\documentclass[a4paper,12pt]{article}
\usepackage{amsmath,amsfonts,amssymb,amsthm,amscd}
\usepackage[arrow, matrix, curve]{xy}
\usepackage[colorlinks=true,citecolor=blue]{hyperref}

\font\cmssl=cmss10 at 12 pt

\usepackage{slashed}
\usepackage[normalem]{ulem}

\newtheorem{thm}{Theorem}
\newtheorem{lem}[thm]{Lemma}
\newtheorem{prop}[thm]{Proposition}
\newtheorem{defn}[thm]{Definition}
\newtheorem{cor}[thm]{Corollary}
\newtheorem{rem}[thm]{Remark}
\newtheorem{notation}[thm]{Notation}
\newtheorem{exa}[thm]{Example}

\title{Fiber-wise linear $F$-manifolds,  compatible  flat connections and   Euler fields}

\date{\today}

\author{Liana David}

\begin{document}
\maketitle

\begin{abstract}
We define  a fiber-wise linear (shortly, FWL)  $F$-manifold as an $F$-manifold
on the total space of a vector bundle 
$\pi : E \rightarrow M$ for which the multiplication and Euler field are FWL 
tensor fields. 
We develop a duality between  FWL  $F$-manifolds 
on $E$ and the total space  $E^{*}$ of the dual vector bundle and we enrich it with FWL Euler fields and FWL compatible connections. 
We present 
examples in dimensions two and three. 
We construct prolongations of $F$-manifolds and  we prove that if the initial $F$-manifold admits a suitable Euler field or compatible flat connection, then all  prolongations inherit FWL Euler fields and  FWL compatible flat connections. 
\end{abstract}

{\it MS  Classification}: 53D45,  53C15.

\section{Introduction}

$F$-manifolds were defined by Hertling and Manin  \cite{h-m}.
The next definition collects  the  notions from the theory of $F$-manifolds to be used in this paper.

\begin{defn} i)  Let $(M, \circ , e)$ be a  (real or complex) manifold together with a 
tensor field  $\circ$ of type $(2,1)$ on $M$, which defines a 
commutative, associative multiplication on $TM$  with  unit field $e$. Then $(M, \circ , e )$ is an  {\cmssl $F$-manifold} if, for any vector fields $X$, $Y$, 
\begin{equation}\label{integr-F-man}
\mathcal L_{X\circ Y} (\circ ) = X\circ  {\mathcal L}_{Y} (\circ ) + Y\circ \mathcal L_{X} (\circ ),
\end{equation}
where $\mathcal L_{X}$ denotes the Lie derivative in the direction of  $X$.\

ii) A vector field $\mathcal E$ is called an {\cmssl Euler field} on an $F$-manifold $(M, \circ , e)$ if  $\mathcal L_{\mathcal E} (\circ ) = \circ$. That is,  for any vector fields $X$, $Y$,
$$
\mathcal L_{\mathcal E} (X\circ Y) -  (\mathcal L_{\mathcal E} X)\circ Y -  X\circ (\mathcal L_{\mathcal E} Y) = X\circ Y.
$$
iii) A torsion-free connection $\nabla$ on $M$ is called {\cmssl compatible} on an $F$-manifold $(M, \circ , e)$ if  $\nabla e =0$ and, for any vector fields $X, Y, Z$, 
$$
\nabla_{X} (\circ )(Y, Z) = \nabla_{Y} (\circ ) (X, Z).
$$ 
iv) A  {\cmssl flat $F$-manifold}  is an $F$-manifold $(M, \circ, e)$ together with a flat compatible connection $\nabla$.
A {\cmssl flat $F$-manifold with Euler field} is a flat $F$-manifold $(M, \circ , e, \nabla )$ together with an Euler field ${\mathcal E}$
such that   $\nabla^{2} {\mathcal E} =0$, where, for any vector fields $X$, $Y$, 
\begin{equation}
(\nabla^{2} {\mathcal E} )_{X, Y} := \nabla_{X} \nabla_{Y} \mathcal E - \nabla_{\nabla_{X} Y} \mathcal E  .
\end{equation}
\end{defn}

$F$-manifolds represent  an important topic with ramifications in various areas of mathematics. 
Any Frobenius manifold without metric is an $F$-manifold \cite{dub}. For an overview on the relation between Frobenius and $F$-manifolds and 
meromorphic connections  see e.g. \cite{h-book, sabbah}. 
Generically semisimple $F$-manifolds with smooth spectral cover  are completely 
understood from the point of view of isolated hypersurface singularities \cite[Theorems 5.3 and 5.6]{h-book}. The integrability condition 
(\ref{integr-F-man})  can be expressed in terms of the  standard Poisson structure on $T^{*}M$, making a link between the theory of $F$-manifolds 
and symplectic geometry \cite[Theorem 2.5]{h-t}. Flat $F$-manifolds with Euler fields were constructed on orbit spaces 
of finite complex  reflection groups  in  \cite{AL17,KMS15,KMS18} and found applications in integrable systems. 

With the motivation to construct new examples, 
in this paper we  define and investigate the  class of  fiber-wise  linear (shortly, FWL) 
$F$-manifolds.
Recall that a tensor field $T\in \mathcal T^{p,q} (E)$ on the total space of a vector bundle 
$\pi : E \rightarrow M$ 
is  FWL
if  $h_{t}^{*} T = t^{1-q} T$
for any $t\neq 0$, 
where $h_{t}: E \rightarrow E$ denotes the  multiplication by the scalar $t.$ 
A connection $\nabla$ on the manifold  $E$ is FWL if $\nabla_{X}Y$ is  a FWL vector field, whenever $X$ and $Y$ are FWL vector  fields \cite{pul}.

 \begin{defn} A  {\cmssl FWL $F$-manifold}
 is an $F$-manifold  
 $(E, \circ , e)$ for which $\circ$ and $e$ are FWL tensor fields.
 A {\cmssl FWL $F$-manifold with Euler field (or compatible connection)} is a FWL $F$-manifold  together with a FWL Euler  field
 (or FWL compatible connection, respectively). 
 \end{defn}

{\bf Outline of the paper.} In Section \ref{preliminary-sect}, intended to fix notation, 
we review 
basic facts on FWL tensor fields and connections.
In  Proposition \ref{comp-prop}, which is a counterpart for  the description of  FWL vector-valued  $2$-forms using generalized derivations  of degree two \cite{b-d},  we  review the description of   FWL symmetric $(2,1)$-tensor fields in terms of  their  components  $(D, l, T^{M})$, where 
$D: \Gamma (E) \rightarrow \Gamma (S^{2} T^{*}M)$, $l : E \rightarrow T^{*}M\otimes E$ and $T^{M}\in \Gamma ( S^{2} T^{*}M\otimes
TM)$ satisfy  a Leibniz type relation \cite{pul}. We recall, following \cite{pul},  the description of  FWL  torsion-free  connections $\nabla$
on $E$ in terms of  data $(\nabla^{M}, \nabla^{E}, \mathcal S )$ (called the components of $\nabla$) where $\nabla^{M}$  (the basic component) 
is a torsion-free connection
on $M$, $\nabla^{E}$ is a connection on  the vector bundle $E$ and $\mathcal S\in  \Gamma (  S^{2}T^{*}M\otimes \mathrm{End}\, E)$.
In Subsection \ref{FWL-sect}  we  prove various facts  on  FWL torsion-free connections, to be used  in the next sections.

In Section \ref{components-section}  we 
develop a general theory for  FWL $F$-manifolds in terms of components of the structures involved. 
In Subsection \ref{components-Fman}
we consider a FWL symmetric  tensor field $\circ$ of type $(2,1)$ and a FWL vector field $e$ on $E$  and we 
determine general conditions such that $(E, \circ , e)$ is a (FWL) $F$-manifold
(see Lemmas \ref{lem-asoc}, \ref{lem-unit} and \ref{prop-integr}).  In particular, we deduce  that any
FWL $F$-manifold $(E, \circ , e)$ lies over an $F$-manifold $(M, *, \bar{e})$, where 
$*$ is the basic component of $\circ$ and 
 the unit field $\bar{e}\in {\mathfrak X}(M)$ 
is the projection of the unit field $e\in {\mathfrak X}_{\mathrm{fwl}} (E)$.
In Subsection  \ref{euler-sect} we consider a FWL vector field $\mathcal E\in {\mathfrak X}_{\mathrm{fwl}} (E)$ on a
FWL $F$-manifold $(E, \circ , e)$ and we determine  
general  conditions  such that 
$\mathcal E$ is an Euler field on $(E, \circ , e)$ (see Proposition \ref{euler-cor}).  In particular, we deduce that 
any FWL Euler field  $\mathcal E$ on $(E, \circ , e)$  projects to an Euler
field $\bar{\mathcal E}$ on $(M, *, \bar{e}).$   In Subsection \ref{conn-sect} we determine 
general conditions on  a FWL torsion-free connection $\nabla$  on $E$ to be  compatible on $(E, \circ , e)$
(see Propositions \ref{nabla-circ-sym}, \ref{mathcal-e} and Lemma \ref{mathcal-s}).  
In particular, we deduce that the basic component $\nabla^{M}$  of $\nabla$  is  compatible on $(M, *, \bar{e}).$\

In Section \ref{duality-section} we 
fix a connection $\nabla$ on $M$ and we use it to define a duality 
between FWL $F$-manifolds on $E$ and  $E^{*}$ 
(see Theorem \ref{short} below). 
Considering the symmetric derivative  
defined by  $\nabla$ on symmetric tensor fields
 as the analogue of the exterior derivative on differential forms,  our duality is analogous to the isomorphism  between the spaces of generalized derivations  on $E$ and $E^{*}$ \cite{d-n}.
More precisely, given a FWL  symmetric $(2,1)$-tensor field 
$\circ$  on $E$, we define, in terms of components, using $\nabla$, 
a FWL symmetric $(2,1)$-tensor field  $\bar{\circ}$ on $E^{*}$ 
(see the beginning of Section \ref{duality-section})  and we
 determine general conditions such that 
\begin{equation}\label{map}
(E, \circ , e)\mapsto (E^{*}, \bar{\circ}, e^{*})
\end{equation}
maps $F$-manifolds to $F$-manifolds  (see Lemmas \ref{duality-asoc}, \ref{duality-unit} and \ref{duality-F}).
The   resulting conditions  
are  involved in full generality,
but they are always satisfied when
$\nabla$ makes  $(M, *, \bar{e})$ a flat $F$-manifold. We obtain (see also Definition \ref{over}):  
\begin{thm}\label{short} Let  
$\pi : E \rightarrow M$ be a vector bundle with base a flat $F$-manifold
$(M, *, \bar{e}, \nabla )$.  There is a bijection   
(induced by the map (\ref{map}))
\begin{equation}\label{i-nabla}
I^{\nabla} : \mathcal F_{\mathrm{fwl}} (E, *, \bar{e} ) \rightarrow  \mathcal F_{\mathrm{fwl}} (E^{*}, *, \bar{e} )
\end{equation}
between   the sets  $ \mathcal F_{\mathrm{fwl}} (E, *, \bar{e} ) $ and 
$\mathcal F_{\mathrm{fwl}} (E^{*}, *, \bar{e} )$ of FWL $F$-manifolds
on $E$ and $E^{*}$ respectively, over   $(M, *, \bar{e})$.
\end{thm}

In Section \ref{duality-enrich} we enrich Theorem 
\ref{short}  with Euler fields and compatible connections. 
In Subsection \ref{duality-euler-sect}
we determine general conditions such that 
the dual $\mathcal E^{*}$ 
(see Remark \ref{various-basic} ii))
of an Euler field on $(E, \circ , e, \mathcal E)$ is an Euler field on 
$I^{\nabla}(E, \circ , e)$ (see Proposition \ref{euler-duality}). 
In Subsection \ref{duality-conn-sect}    we  show that any FWL torsion-free connection  
$\nabla^{fwl}$ 
on $E$ defines naturally a FWL torsion-free connection 
$(\nabla^{fwl})^{*}$ on $E^{*}$  (see Definition \ref{dual-connection}) and  we compare the properties of $\nabla^{fwl}$ and $(\nabla^{fwl})^{*}$ 
(see Proposition \ref{compat-dual}  and Lemma \ref{twice}).
Combining (\ref{map}) with  the maps $\mathcal E\mapsto \mathcal E^{*}$ and $\nabla^{fwl}\mapsto  (\nabla^{fwl})^{*}$ we obtain 
(see Definitions \ref{over-euler} and  \ref{over-conn}):

\begin{thm}\label{duality-euler-thm}
Let  $\pi : E \rightarrow M$ be a vector bundle with base a flat $F$-manifold  $(M, *, \bar{e}, \nabla ).$\

i) There is a bijection  between the sets of  FWL  flat $F$-manifolds  
on $E$ and $E^{*}$  over  $(M, *, \bar{e}, \nabla )$.\

ii) Assume that $(M, *, \bar{e}, \nabla ,\bar{\mathcal E})$ is a flat $F$-manifold with Euler field. 
There is a bijection   between the sets of FWL $F$-manifolds with Euler fields 
(respectively, FWL flat $F$-manifolds with Euler fields)  on $E$ and $E^{*}$ over  $(M, *, \bar{e}, \bar{\mathcal E})$
(respectively, over $(M, *, \bar{e}, \nabla , \bar{\mathcal E})$). 
\end{thm}

In Section \ref{examples-sect}  we present examples in small dimensions,  inspired by the local classification of complex 
$F$-manifolds in dimension two and three  \cite{h-paper,h-book}.   
In Section \ref{tangent-sect} we define and investigate  the tangent prolongation of an $F$-manifold. 
We start  by  proving:

\begin{prop}\label{integr-comp}
Let $(M, * , e ) $ be an $F$-manifold. The data $(D, l, r := * )$ where 
\begin{equation}\label{prolong} 
(DX) (Y, Z) := \mathcal L_{X} ( *  ) (Y, Z),\ ( lX) (Y) := X *  Y,\ r(X, Y) := X*Y
\end{equation}
are the components of a FWL multiplication $ *^{T} $ on $TM$, which,  together with the tangent prolongation  $e^{T}$ of $e$
(see Remark \ref{various-basic} iii)),
defines a FWL $F$-manifold $(TM,  *^{T} , e^{T}).$  We call it the {\cmssl tangent prolongation}
of $(M,  *  , e).$
\end{prop}

Then we add to Proposition \ref{integr-comp}  Euler fields and compatible connections.
Recall first 
that any  torsion-free  connection $\nabla$ on $M$ induces a canonical FWL torsion-free connection $\nabla^{T}$ on $TM$ 
(called the tangent prolongation or the complete lift of $\nabla$, in the language of \cite{t-c}),
which satisfies $\nabla^{T}_{X^{T}} (Y^{T}) = (\nabla_{X} Y)^{T}$
for any $X, Y\in {\mathfrak X}(M).$  
After determining  the components of  $\nabla^{T}$ (see Lemma \ref{complete}) we investigate when
$\nabla^{T}$ is compatible on $(TM, *^{T}, e^{T})$ (see Proposition \ref{compat-tg}).  As a consequence, we obtain:

\begin{prop}\label{e-con}  i) 
If $\mathcal E$ is an Euler field on an $F$-manifold $(M, *, e)$, then  its tangent prolongation $\mathcal E^{T}$ is an Euler field on $(TM,  *^{T} , e^{T}).$\
 
ii) If $(M, *, e, \nabla )$ is a flat  $F$-manifold,  then  so is  $(TM,  *^{T} , e^{T}, \nabla^{T}) $.\

iii)  If $(M, *, e, \nabla  , \mathcal E)$ is a flat  $F$-manifold  with Euler field, then  so is   $(TM,  *^{T} , e^{T},\nabla^{T},  \mathcal E^{T})$.
\end{prop}
In Section \ref{ctg-def} we define the 
cotangent and generalized prolongation of a flat $F$-manifold $(M, *, e, \nabla )$. 
The latter  makes the generalized tangent bundle $\mathbb{T}M = TM\oplus T^{*}M$  an $F$-manifold and provides a link between generalized geometry and the theory of $F$-manifolds.   This view-point  will be investigated more closely  in another work. 
We explain how the tangent, cotangent and generalized prolongations  can be iterated, providing 
further examples of FWL $F$-manifolds, with  Euler fields and flat compatible connections  inherited from the initial $F$-manifold
(see Corollary  \ref{last}).

\section{Preliminary material}\label{preliminary-sect}

\begin{notation}{\rm   
All results are stated in the smooth category.
We denote by ${\mathfrak X}(M)$,  $\Omega^{k}(M)$  and $\mathcal T^{p, q}(M)$ the sheaves of vector fields,  $k$-forms and tensor fields of type $(p, q)$ on a manifold $M$.  Given a connection $\nabla$ on $M$ we define, following \cite{rubio}, 
the symmetric bracket  $[\cdot , \cdot ]_{\nabla} :{\mathfrak X}(M) \times  \mathfrak{X}(M)\rightarrow {\mathfrak X}(M)$ by 
$[X, Y]_{\nabla}:= \nabla_{X}Y +\nabla_{Y} X$, for any   $X, Y\in {\mathfrak X}(M)$.
For  $X\in {\mathfrak X} (M)$,  the Lie derivative $\mathcal L_{X}\nabla\in \mathcal T^{2,1}(M)$  is defined by
$$
(\mathcal L_{X}\nabla )_{Y} Z:= \mathcal  L_{X} (\nabla_{Y} Z) - \nabla_{\mathcal L_{X} Y} (Z) -\nabla_{Y} ( \mathcal L_{X} Z),\ 
\forall  Y, Z\in {\mathfrak X}(M).$$
We denote by   ${\mathfrak X}_{\nabla}(M)$  the sheaf of
$\nabla$-parallel vector fields with respect to a flat connection $\nabla .$ 
If $\pi : E\rightarrow M$ is a vector bundle and $X\in {\mathfrak X}(E)$ is projectable, 
we often denote by $\bar{X}\in {\mathfrak X}(E)$ its projection to $M$.
We denote by $\Omega^{k}(M,E)$ the space of $k$-forms on $M$ with values in $E$.}
\end{notation} 

\subsection{FWL tensor fields}

Let  $\pi : E\rightarrow M$ be a vector bundle, $(x^{i})$  local coordinates on $M$, $( s_{i} )$ a local trivialization of $E$ and $(x^{i}, \xi^{j})$ the induced local coordinates on $E$.
Locally, FWL  vector fields are $C^{\infty}(M)$ linear combinations of vector fields of the form
$\frac{\partial}{\partial x^{i}}$ and $\xi^{j} \frac{\partial}{\partial \xi^{ i}}$,  FWL  $1$-forms are $C^{\infty}(M)$ linear combinations of  
$1$-forms of the form 
$d\xi^{i}$ and $\xi^{j} dx^{i}$  and  FWL $(2,1)$-tensor fields are of the form
\begin{align}
\nonumber& T = \left( a^{i}_{jk} (x) d\xi^{j}\otimes dx^{k} + \tilde{a}^{i}_{kj} (x) dx^{k} \otimes d\xi^{j}  +\xi^{j} a^{i}_{j, kp} (x) dx^{k} \otimes dx^{p} \right) 
\otimes \partial  \xi^{i} \\
\label{loc}&  + b^{a}_{ij}(x) dx^{i}\otimes dx^{j} \otimes \partial x^{a},
\end{align}
where  $a^{i}_{jk}, \tilde{a}^{i}_{kj}, a^{i}_{j, kp}, b_{a, ij}\in C^{\infty}(M)$ and to simplify 
notation $\frac{\partial}{\partial x^{a}}$ is denoted by $\partial x^{a}$ (similarly for 
$\frac{\partial}{\partial \xi^{i}}$).
The $C^{\infty }(M)$-modules of FWL, FWL   symmetric (in the first $p$-arguments) tensor fields of type $(p,q)$ 
and FWL vector fields 
are  denoted by $\mathcal T^{p,q}_{\mathrm{fwl}}(E)$, $\mathcal T^{p,q}_{\mathrm{fwl},s}(E)$
and  $\mathfrak{X}_{\mathrm{fwl}} (E)$ respectively.\

A tensor field $T\in {\mathcal T}^{p,q}(E)$ is   {\cmssl core}  if $h_{t}^{*} T = t^{-q} T$, for any $t\neq 0$, where $h_{t}$ is the multiplication by the scalar $t$ 
(see \cite[Section 5] {pul} for details on core vector fields).
We denote by
$\mathcal T^{p,q}_{\mathrm{core}} (E)$
and ${\mathfrak X}_{\mathrm{core}} (E)$ 
the  $C^{\infty}(M)$-modules of core  tensor fields of type $(p,q)$ and core vector fields respectively.  A  $p$-form is core if  and only if it is the pullback of a $p$-form on $M$. 
The vertical lift $s^{\uparrow }\in {\mathfrak X}(E)$ of a section $s\in \Gamma (E)$, defined 
by $s^{\uparrow}(p) = \frac{d}{dt}\vert_{t=0} ( p+  t s(\pi p))$ ($p\in E$), is a core vector field and any core vector field is of this form.
This defines a  $C^{\infty}(M)$-linear isomorphism ${\mathfrak X}_{\mathrm{core}} (E)\cong \Gamma (E)$ which extends to an isomorphism 
 \begin{equation}\label{iso-core}
 \mathcal T^{p,1}_{\mathrm{core}} (E) \cong \Gamma ( (T^{*}M)^{p}\otimes  E)
 \end{equation}
 of $C^{\infty}(M)$ modules ($p\in \mathbb{N}$). Any   $T\in {\mathcal T}^{p,1}_{\mathrm{core} }(E)$   is  locally of the form
$$
 T = T^{k}_{i_{1}\cdots i_{p}} (x)  dx^{i_{1}} \otimes \cdots \otimes dx^{i_{p}} \otimes \partial \xi^{k},
 $$
 where $T^{k}_{i_{1}\cdots i_{p}} \in C^{\infty}(M)$.  For any 
 $\omega \in \Gamma ( (T^{*}M)^{p}\otimes  E)$, we denote by $\omega^{\uparrow}\in \mathcal T^{p,1}_{\mathrm{core}} (E)$   
 the corresponding core tensor field in the isomorphism (\ref{iso-core}).\

\begin{lem}\label{lin-core} \cite{mckanzie,pul} i) A tensor field $T\in {\mathcal T}^{p,1} (E)$ is FWL if and only if $T(X_{1}, \cdots , X_{p} ) \in {\mathfrak X}_{\mathrm{fwl}} (E)$,
for any $X_{i} \in  {\mathfrak X}_{\mathrm{fwl}} (E)$.\

ii)  Let $T\in {\mathcal T}^{p,1}_{\mathrm{fwl}} (E)$ and $s\in \Gamma (E).$ 
The natural contractions  between $T$ and  $s^{\uparrow}$ ($p$ in total) 
and the Lie derivative $\mathcal L_{s^{\uparrow}} T$ are core tensor fields.\

iii) Any  FWL vector field is projectable and the map $\pi_{*} : {\mathfrak X}_{\mathrm{fwl}} (E) \rightarrow  {\mathfrak X} (M)$ is surjective.  The Lie bracket of FWL vector fields is a FWL vector field. 
\end{lem}

\begin{rem}\label{various-basic}{\rm   i) A derivation  is a map 
 $\Delta : \Gamma (E) \rightarrow \Gamma (E)$ which satisfies
 $$
 \Delta (fs) = f \Delta (s) +  a(\Delta)  (f) s,\ \forall  f\in C^{\infty}(M),\ s\in \Gamma (E),
 $$
 where $a(\Delta )\in {\mathfrak X}(M)$ is   the symbol of $\Delta$.
 If $X\in{ \mathfrak X}_{\mathrm{fwl}} (E)$ and $s\in \Gamma (E)$, then $\mathcal L_{X} s^{\uparrow}\in {\mathfrak X}_{\mathrm{core}} (E)$ 
 (from Lemma \ref{lin-core} ii))
 and 
 \begin{equation}\label{deriv-vector-field}
 \mathcal L_{X} s^{\uparrow} =  (\Delta_{X} s)^{\uparrow},
  \end{equation}
 where $\Delta_{X}$ is a derivation with   $a(\Delta_{X}) = \pi_{*} X$.
 The map $X\mapsto \Delta_{X}$ is a Lie algebra isomorphism  between  FWL vector fields and derivations on $E$
 \cite{mckanzie}:
  \begin{equation}\label{comutator}
 \Delta_{\mathcal L_{X}Y} = [\Delta_{X}, \Delta_{Y}],\ \forall X, Y\in {\mathfrak X}_{\mathrm{fwl}} (E). 
 \end{equation}
ii) There is a Lie algebra isomorphism  ${\mathfrak X}_{\mathrm{fwl}} (E)\ni X\mapsto X^{*}\in {\mathfrak X}_{\mathrm{fwl}} (E^{*})$, where
$$
(\Delta_{X^{*}} \mu ) s := ( \pi_{*} X)  (\mu (s))  - \mu (\Delta_{X} s),\ \forall \mu \in \Gamma (E^{*}),\ s\in \Gamma (E).
$$ 
iii) For any  $X\in {\mathfrak X}(M)$,   $\mathcal L_{X}\in \mathrm{End}\,  {\mathfrak X}(M)$ is a derivation on $TM$ with symbol $X$.  It defines the {\cmssl tangent prolongation $X^{T} \in {\mathfrak X}_{\mathrm{fwl}}(TM)$ of $X$.}} \end{rem}

 \begin{prop}\label{comp-prop} A 
 FWL symmetric $(2,1)$-tensor field $T$  is uniquely determined by the data 
 $(D, l ,T^{M})$ (called {\cmssl components of $T$}) where 
 \begin{equation}\label{components}
D :\Gamma (E) \rightarrow \Gamma (  S^{2} T^{*}M\otimes E),\   l: E \rightarrow T^{*}M\otimes E
\end{equation}
 and $T^{M} \in \Gamma (  S^{2} T^{*}M\otimes TM)$ (called the {\cmssl basic component of $T$}) 
are given by 
\begin{equation}\label{data-1}
(D s)^{\uparrow} := \mathcal L_{s^{\uparrow}} T,\  (l s)^{\uparrow} :=  T(s^{\uparrow}, \cdot ),\ 
T^{M} (X_{1},  X_{2}) := \pi_{*}  T (\tilde{X}_{1},  \tilde{X}_{2} )
\end{equation}
where $\tilde{X}_{i}\in  {\mathfrak X}(E)$ project to $X_{i} \in {\mathfrak X}(M)$. 
The  Leibniz-type relation  
  \begin{equation}\label{comp-cond}
 D_{X, Y} (fs) = f D_{X, Y} s +  X(f) l_{Y} s+ Y(f) l_{X}s  -  df ( T^{M}(X, Y) ) s,
 \end{equation}
 is satisfied, where $f\in C^{\infty}(M)$,  $s\in \Gamma (E)$ and $X, Y\in {\mathfrak X}(M)$.  
 Conversely, any such data corresponds to a unique 
 FWL symmetric $(2,1)$-tensor field.
  \end{prop}

\begin{proof} 
Let $T\in {\mathcal T}^{2,1}_{\mathrm{fwl}}(E)$.
Using Lemma \ref{lin-core} ii) and the isomorphism (\ref{iso-core}), we define $D$ and $l^{(i)}$ 
($1\leq i\leq 2$)  by
\begin{equation}\label{data-1}
(D s)^{\uparrow} := \mathcal L_{s^{\uparrow}} T,\  (l^{(1)}s)^{\uparrow} :=  T(s^{\uparrow}, \cdot ),\ (l^{(2)}s)^{\uparrow} :=  T(\cdot  , s^{\uparrow})
\end{equation}
and  $T^{M}$ by
\begin{equation}\label{data-2}
T^{M} (X_{1},  X_{2}) := \pi_{*}  T (\tilde{X}_{1},  \tilde{X}_{2} )
\end{equation}
where $\tilde{X}_{i}\in  {\mathfrak X}(E)$ project to $X_{i} \in {\mathfrak X}(M)$. 
If $T$ is given  by  (\ref{loc}), then 
\begin{align}
\nonumber& Ds_{j} = a^{i}_{j, kp}dx^{k} \otimes dx^{p} \otimes s_{i},\   l^{(1)} s_{j}= a^{i}_{jk} dx^{k}\otimes s_{i},\\
\label{loc-1}& l^{(2)} s_{j} = \tilde{a}^{i}_{kj} dx^{k}\otimes s_{i},\  T^{M} = b^{a}_{ij}dx^{i}\otimes dx^{j} \otimes \frac{\partial}{\partial x^{a}}
\end{align}
and 
\begin{equation}
a^{i}_{jk} = s_{i}^{*} ( l^{(1)}_{ \partial  x^{k}} s_{j}),\ \tilde{a}^{i}_{kj} = s_{i}^{*} ( l^{(2)}_{\partial x^{ k}} s_{j}),
a^{i}_{j, kp} := s_{i}^{*} ( D_{\partial x^{ k}, \partial x^{ p}} s_{j}),
\end{equation}
where $(s_{i}^{*})$ is the dual trivialization of $( s_{i})$ and  to simplify notation we write $\partial x^{ k}$ instead of $\frac{\partial}{\partial x^{k}}$. 
Now, $T$ is symmetric  if and only if
$a^{i}_{j, kp} = a^{i}_{j, pk}$ 
(i.e. $D$ takes values in $\Gamma (S^{2} T^{*}M\otimes E)$),  $T^{M}$ is symmetric 
and $l^{(1)} = l^{(2)}.$ Let $l:= l^{(1)} = l^{(2)}$.
Relation (\ref{comp-cond})  follows from the  definitions of $D$, $l$, $T^{M}$ and  
\begin{align}
\nonumber& (\mathcal L_{fZ}T )(X_{1},  X_{2}) = f (\mathcal L_{Z}T )(X_{1}, X_{2})\\
\nonumber&  + X_{1} (f) T (Z, X_{2})  + X_{2}(f) T (X_{1}, Z)  -  df (T (X_{1}, X_{2} ))Z
\end{align}
for any $X_{1}, X_{2}, Z\in {\mathfrak X}(E)$ and $f \in C^{\infty}(E).$
\end{proof}

\begin{rem}\label{iso-sum}{\rm  Let  $\pi_{i} : E_{i} \rightarrow M$ be two vector bundles over the same manifold.\
If $X_{i} \in {\mathfrak X}_{\mathrm{fwl}}(E_{i})$ have the same projection on $M$, then  their direct sum 
 $X_{1}\oplus X_{2} \in {\mathfrak X}_{\mathrm{fwl}} (E_{1}\oplus E_{2})$ is defined by the derivation
 $\Delta_{X_{1}\oplus X_{2}} (s_{1}, s_{2}) : = (\Delta_{X_{1}} s_{1}, \Delta_{X_{2}} s_{2})$.
 If $T_{i} \in \mathcal T^{2,1}_{\mathrm{fwl},s} (E_{i})$ have components  $(D^{i}, l^{i}, T^{M})$  then their direct sum
$T_{1} \oplus T_{2} \in \mathcal T^{2,1}_{\mathrm{fwl},s} (E_{1}\oplus E_{2})$ has
components $(D^{1}\oplus D^{2}, l^{1}\oplus l^{2}, T^{M})$, where, for any $s_{i}\in \Gamma (E_{i})$ and $X, Y\in {\mathfrak X}(M)$, 
\begin{align*}
\nonumber& (D^{1}\oplus D^{2})_{X, Y }(s_{1}, s_{2}) = (D^{1}_{X, Y }(s_{1}), D^{2}_{X, Y }(s_{2})),\\
\nonumber& (l^{1}\oplus l^{2})_{X}(s_{1}, s_{2}) = (l^{1}_{X}(s_{1}), l^{2}_{X}(s_{2})).
\end{align*}}
\end{rem}

\begin{rem}{\rm  Proposition \ref{comp-prop} can be extended to FWL
(not necessarily symmetric) 
tensor fields of type $(p,1)$ for any $p\in \mathbb{Z}_{\geq 1}.$ 
Given $T\in \mathcal T^{p,1}_{\mathrm{fwl}} (E)$  one defines as in   (\ref{data-1}) and (\ref{data-2})
$(D, l^{(i)}, T^{M})$ (called components of $T$), 
where $D: \Gamma (E) \rightarrow \Gamma ( (T^{*}M)^{p}\otimes E)$, $l^{(i)}: E \rightarrow (T^{*}M)^{p-1} \otimes E$ 
($1\leq i\leq p$)  and
$T^{M}\in \Gamma ( (T^{*}M)^{p}\otimes TM)$  satisfy a Leibniz type relation  similar to (\ref{comp-cond}) and determine $T$ uniquely. 
If $M$ has a flat connection $\nabla$, then $T =0$ if and only if $D_{{X}_{1}, \cdots , X_{p}} =0$, $l^{(i)}_{X_{1}, \cdots , X_{p-1}} =0$
(for any $1\leq i\leq p$) 
and $T^{M} (X_{1}, \cdots , X_{p}) =0$ for any $X_{i} \in {\mathfrak X}_{\nabla}(M)$.}
\end{rem}

\subsection{FWL torsion-free connections}\label{FWL-sect}

A  FWL torsion-free connection
 $\nabla$ on $E$ is uniquely determined by the data
 $(\nabla^{M}, \nabla^{E}, \mathcal S )$ 
  (called {\cmssl components of $\nabla$}) where $\nabla^{M}$ is a torsion-free connection $M$,
  $\nabla^{E}$ is a connection on the vector bundle $E$ and $\mathcal S \in \Gamma ( S^{2} T^{*}M\otimes \mathrm{End}\, E)$
  (see Corollary 8.8 of \cite{pul}). Conversely, any  such data corresponds to a unique 
 FWL torsion-free connection. 
 The  basic component $\nabla^{M}$ is defined by
 $\nabla^{M}_{X}Y = \pi_{*} \nabla_{\tilde{X}} \tilde{Y}$, for any $\tilde{X}, \tilde{Y}\in {\mathfrak X}(E)$ which project to $X, Y\in{\mathfrak X}(M)$. 
 To define  $\nabla^{E}$ one notices  that  $\nabla s^{\uparrow}$ is a  $(1,1)$-tensor field  of core type, 
 for any $s\in \Gamma (E)$. Then  $\nabla^{E} s\in \Omega^{1}(M, E)$  is defined by 
 \begin{equation}\label{fwl-E}
(\nabla^{E}_{X}  s)^{\uparrow} =  \nabla_{\tilde{X}}  (s^{\uparrow}).
 \end{equation}
 Finally, $\mathcal L_{s^{\uparrow}} \nabla\in \mathcal T^{2,1}_{\mathrm{core}}(E)$ determines $\mathcal S s\in \Gamma (S^{2} T^{*}M\otimes E)$ 
 by 
 \begin{equation}\label{g-s}
(\mathcal L_{s^{\uparrow}} \nabla )(\tilde{X}, \tilde{Y}) =\left( \frac{1}{2} \left( \nabla^{E}_{X}\nabla^{E}_{Y} +\nabla^{E}_{Y} \nabla^{E}_{X}
 -\nabla^{E}_{[X, Y]_{\nabla^{M}}} \right)  s+  \mathcal S_{X ,Y}s\right)^{\uparrow}.
 \end{equation}
 Let $\mathcal G$ be the expression arising in the right hand side of (\ref{g-s}). Thus,
  \begin{equation}\label{deriv-G}
 (\mathcal L_{s^{\uparrow}}\nabla )(\tilde{X}, \tilde{Y})
 =( \mathcal G_{X, Y}s)^{\uparrow}
 \end{equation}
 and 
 \begin{align}
\nonumber \mathcal G_{X, Y} s &= \frac{1}{2} \left( \nabla^{E}_{X}\nabla^{E}_{Y} +\nabla^{E}_{Y} \nabla^{E}_{X}
 -\nabla^{E}_{[X,Y]_{\nabla^{M}}} \right)  s+  \mathcal S_{X ,Y}s\\
 \label{g-s-1}& = (\nabla^{E})^{2}_{X, Y} s -\frac{1}{2} R^{\nabla^{E}}_{X, Y} s +  \mathcal S_{X ,Y }s
\end{align}
where $(\nabla^{E})^{2}_{X, Y} s := \nabla^{E}_{X} \nabla^{E}_{Y }s - \nabla^{E}_{\nabla^{M}_{X} Y}s$.
We often refer to $(\nabla^{M}, \nabla^{E},\mathcal G  )$ (instead of $(\nabla^{M} , \nabla^{E}, \mathcal S  )$) as  components 
 of $\nabla.$ 
 
 \begin{lem}\label{prima-nabla} Let $X, Y\in {\mathfrak X}_{\mathrm{fwl}} (E)$ and  $s\in \Gamma (E)$. Then 
 \begin{align}
\nonumber &\nabla_{s^{\uparrow}} X = (\nabla^{E}_{\bar{X}} s - \Delta_{X} s)^{\uparrow}\\
\nonumber& \Delta_{\nabla_{X} Y} = - \mathcal G_{\bar{X}, \bar{Y}} s +\nabla^{E}_{\bar{X}} \Delta_{Y} s  + \nabla^{E}_{\bar{Y}} \Delta_{X} s 
 - \Delta_{Y} \Delta_{X} s.\\
\label{need-nabla}& [ \nabla_{\bar{X}}^{E}, \Delta_{\mathcal E} ] s= \nabla^{E}_{\mathcal L_{\bar{X}} \bar{\mathcal E}} s.
 \end{align}
 \end{lem}
 
 \begin{proof}  The first relation in (\ref{need-nabla}) uses  the torsion-free property of $\nabla$:
 \begin{equation}
 \nabla_{s^{\uparrow}} X = \nabla _{X} (s^{\uparrow} ) + \mathcal L_{s^{\uparrow}} X = (\nabla^{E}_{\bar{X}} s)^{\uparrow} -
( \Delta_{X} s)^{\uparrow}.
 \end{equation}
 The second relation in (\ref{need-nabla}) follows by writing 
 \begin{align}
\nonumber&  (\Delta_{\nabla_{X} Y} s)^{\uparrow} = \mathcal L_{\nabla_{X} Y} (s^{\uparrow}) = - \mathcal L_{s^{\uparrow}} (\nabla_{X}Y)
 =  - \mathcal L_{s^{\uparrow}} (\nabla )_{X, Y} -\nabla_{\mathcal L_{s^{\uparrow}} X} Z - \nabla_{X}
  \mathcal L_{s^{\uparrow}} Y\\
 \nonumber&  = - (\mathcal G_{\bar{X}, \bar{Y}} s)^{\uparrow} +\nabla_{(\Delta_{X} s)^{\uparrow}} Y +\nabla_{X} (\Delta_{Y} s)^{\uparrow},
 \end{align}
 and using  the first relation in (\ref{need-nabla}) and (\ref{fwl-E}). For the third relation in (\ref{need-nabla}),
 \begin{align}
 \nonumber&  \left(\nabla^{E}_{\bar{X}} \Delta_{\mathcal E} s - \Delta_{\mathcal E} \nabla^{E}_{\bar{X}} s\right)^{\uparrow} 
 =\mathcal L_{X} (\Delta_{\mathcal E} s)^{\uparrow} - \mathcal L_{\mathcal E} (\nabla^{E}_{\bar{X}} s)^{\uparrow} \\
 \nonumber&  = \mathcal L_{X} \mathcal L_{\mathcal E} (s^{\uparrow} ) -\mathcal L_{\mathcal E} \mathcal L_{X} ( s^{\uparrow}) =
 \mathcal L_{\mathcal L_{X} \mathcal E} (s^{\uparrow}).
 \end{align}
 On the other hand, 
 $(\nabla^{E}_{\mathcal L_{\bar{X}} \bar{\mathcal E}} s)^{\uparrow} = \mathcal L_{\mathcal L_{X} \mathcal E} (s^{\uparrow})$ because 
 $\mathcal L_{X} \mathcal E$ is FWL and projects to $\mathcal L_{\bar{X}} \bar{\mathcal E}.$ The third relation in (\ref{need-nabla}) follows.
  \end{proof}
 
 \begin{lem} Let  
 $\mathcal E\in {\mathfrak X}_{\mathrm{fwl}} (E).$  Then $\nabla^{2} \mathcal E\in \mathcal T^{2,1}_{\mathrm{fwl}} (E)$ has  basic component
 $(\nabla^{M})^{2} \bar{\mathcal E}$ and the other components 
 $(D,  l^{(1)}, l^{(2)})$ are given by 
 \begin{equation}\label{comp-nabla-e}
 l^{(1)}_{X} s = \mathcal S_{X, \bar{\mathcal E}} s +\frac{1}{2} R^{\nabla^{E}}_{X, \bar{\mathcal E}} s,\ 
 l^{(2)}_{X} s= R^{\nabla^{E}}_{X, \bar{\mathcal E}} s 
 \end{equation}
 and 
 \begin{align}
 \nonumber D_{X, Y} s &= \nabla_{X}^{E} (
 \mathcal S_{Y, \bar{\mathcal E}}s)-  \nabla_{\bar{\mathcal E}}^{E} ( \mathcal S_{X, Y}s) +
 \mathcal S_{X, \nabla^{M}_{Y} \bar{\mathcal E} } s - \mathcal S_{\nabla^{M}_{X}Y , \bar{\mathcal E}} s 
 +\Delta_{\mathcal E} ( \mathcal S_{X, Y}s)\\
 \nonumber&+\frac{1}{2} ( \nabla_{\bar{\mathcal E}}- \Delta_{\mathcal E} ) R^{\nabla^{E}}_{X, Y} s
 + R^{\nabla^{E}}_{X, \bar{\mathcal E}} (\nabla^{E}_{Y} s) -\nabla^{E}_{(\nabla^{M} )^{2}_{X, Y}\bar{\mathcal E}} s\\
 \label{fwl-E-2}& +\frac{1}{2}\nabla_{X}^{E} ( R^{\nabla^{E}}_{Y, \bar{\mathcal E}} s) -\frac{1}{2} R^{\nabla^{E}}_{X,\nabla^{M}_{Y}  \bar{\mathcal E}} s -\frac{1}{2} R^{\nabla^{E}}_{\nabla^{M}_{X}Y, \bar{\mathcal E}}s,
  \end{align}
 for any $X , Y\in {\mathfrak X}(M)$ and $s\in \Gamma (E).$ 
 \end{lem}
 
 \begin{proof} For any $s\in \Gamma (E)$ and $X\in {\mathfrak X}_{\mathrm{fwl}} (E)$, 
 \begin{align}
 \nonumber& (\nabla^{2}\mathcal E)_{s^{\uparrow}, X}= \nabla_{s^{\uparrow}}\nabla_{X} \mathcal E - 
 \nabla_{\nabla_{s^{\uparrow}}X} \mathcal E
 = (\nabla^{E}_{\nabla^{M}_{\bar{X}}\bar{\mathcal E}}s-\Delta_{\nabla_{X}\mathcal E}s)^{\uparrow}
 -\nabla_{ ( \nabla^{E}_{\bar{X}}s-\Delta_{X}s)^{\uparrow}} \mathcal E
 \end{align}
 where we used that  $\nabla_{X} \mathcal E \in {\mathfrak X}_{\mathrm{fwl}} (E)$ projects to $\nabla^{M}_{\bar{X}} \bar{\mathcal E}$ and 
 the first relation in  (\ref{need-nabla}).  Using the second relation in
 (\ref{need-nabla}) for the term  $\Delta_{\nabla_{X}\mathcal E}s$ we obtain
 \begin{equation}
 (\nabla^{2}\mathcal E )_{s^{\uparrow}, X}  = 
\left( \mathcal S_{\bar{X}, \bar{\mathcal E}} s +\frac{1}{2} R^{\nabla^{E}}_{\bar{X}, \bar{\mathcal E}}s+\nabla^{E}_{\mathcal L_{\bar{X}} \bar{\mathcal E}} s 
 -[\nabla^{E}_{\bar{X}} , \Delta_{\mathcal E} ]s\right)^{\uparrow}
 \end{equation}
 which implies the first relation in  (\ref{comp-nabla-e}),
 as  $[\nabla^{E}_{\bar{X}} , \Delta_{\mathcal E} ]s=\nabla^{E}_{\mathcal L_{\bar{X}} \bar{\mathcal E}}s$
 from the third relation in (\ref{need-nabla}).  The second relation in (\ref{comp-nabla-e}) can be proved similarly.
 It remains to prove  (\ref{fwl-E-2}). For this, let $X, Y\in {\mathfrak X}_{\mathrm{fwl}} (E)$ and write 
 \begin{equation}\label{lie-nabla}
\mathcal L_{s^{\uparrow}} ( \nabla^{2} \mathcal E )_{X, Y}= \mathcal L_{s^{\uparrow}} \left(  (\nabla^{2} \mathcal E)_{X, Y}\right) -   (\nabla^{2}\mathcal E)_{ \mathcal L_{s^{\uparrow}} X, Y}
 -  (\nabla^{2}\mathcal E)_{X,  \mathcal L_{s^{\uparrow}} Y}.
\end{equation}
 Since  $\mathcal L_{s^{\uparrow}} X= -(\Delta_{X} s)^{\uparrow}$ and similarly for $\mathcal L_{s^{\uparrow}} Y$, 
 relations (\ref{comp-nabla-e}) imply 
 \begin{equation}\label{t2-3-nabla}
(\nabla^{2}\mathcal E)_{ \mathcal L_{s^{\uparrow}} X, Y} = 
- ( \mathcal S_{\bar{Y}, \bar{\mathcal E}} \Delta_{X} s +\frac{1}{2} R^{\nabla^{E}}_{\bar{Y}, \bar{\mathcal E}} \Delta_{X} s)^{\uparrow},\  (\nabla^{2}\mathcal E)_{X,  \mathcal L_{s^{\uparrow}} Y} = 
-  (R^{\nabla^{E}}_{\bar{X}, \bar{\mathcal E}} \Delta_{Y} s )^{\uparrow}.
 \end{equation}
 The first term in the right hand side of  (\ref{lie-nabla}) equals 
 \begin{align}
\nonumber&  \mathcal L_{s^{\uparrow}}  \left( \nabla_{X}\nabla_{Y} \mathcal E -\nabla_{\nabla_{X}Y} \mathcal E \right) =
\mathcal L_{s^{\uparrow}}(\nabla )_{X} (\nabla_{Y}\mathcal E) + \nabla_{\mathcal L_{s^{\uparrow}}X} (
 \nabla_{Y}\mathcal E )  +\nabla_{X} \mathcal L_{s^{\uparrow}} (\nabla_{Y} \mathcal E )\\ 
  \nonumber& - \mathcal L_{s^{\uparrow}} (\nabla )_{\nabla_{X} Y}\mathcal E -\nabla_{\mathcal L_{s^\uparrow}(\nabla_{X} Y)}
  \mathcal E -\nabla_{\nabla_{X}Y} (\mathcal L_{s^{\uparrow}}\mathcal E)  .
  \end{align}
  Using (\ref{deriv-G}),  (\ref{need-nabla}), (\ref{deriv-vector-field}) and (\ref{fwl-E}), we obtain that the above expression becomes
  \begin{align}
  \nonumber&( \mathcal G_{\bar{X}, \nabla^{M}_{\bar{Y}} \bar{\mathcal E}} s -\nabla^{E}_{\nabla^{M}_{\bar{Y}}\bar{\mathcal E}} (\Delta_{X}s) 
  -\mathcal G_{\bar{Y}, \bar{\mathcal E}} \Delta_{X} s +\nabla^{E}_{\bar{Y}} \Delta_{\mathcal E} \Delta_{X} s +\nabla^{E}_{\bar{\mathcal E}} \Delta_{Y} \Delta_{X} s -\Delta_{\mathcal E} \Delta_{Y} \Delta_{X} s)^{\uparrow}\\
  \nonumber& + \nabla^{E}_{\bar{X}} \left( \mathcal G_{\bar{Y}, \bar{\mathcal E}} s -\nabla^{E}_{\bar{Y}} \Delta_{\mathcal E} s 
  - \nabla^{E}_{\bar{\mathcal E}} \Delta_{Y} s +\Delta_{\mathcal E} \Delta_{Y} s\right)^{\uparrow}\\
 \nonumber& +  \left( \nabla^{E}_{\nabla^{M}_{\bar{X}} \bar{Y}} \Delta_{\mathcal E} s-\mathcal G_{\nabla^{M}_{\bar{X}}\bar{Y}, \bar{\mathcal E}} s 
 +(\Delta_{\mathcal E} -
 \nabla^{E}_{\bar{\mathcal E}} ) \left( \mathcal G_{\bar{X} , \bar{Y}} s- \nabla^{E}_{\bar{X}} \Delta_{Y} s -\nabla^{E}_{\bar{Y}} \Delta_{X} s
 + \Delta_{Y}\Delta_{X} s\right) \right)^{\uparrow}.
  \end{align}
Combing relations  (\ref{t2-3-nabla})  with the above relation, writing the $\mathcal G$-component  in terms of  the $\mathcal S$-component  (see (\ref{g-s-1})) and using the third relation in (\ref{need-nabla}), we obtain  that  $\mathcal L_{s^{\uparrow}} ( \nabla^{2} \mathcal E )_{X, Y}$ is the core vector field defined  by the section  of $E$ from the right hand side of (\ref{fwl-E-2}) with $X$, $Y$ replaced by $\bar{X}$, $\bar{Y}$. 
 \end{proof}
 
 \begin{cor}\label{cor-E} Let $\mathcal E\in {\mathfrak X}_{\mathrm{fwl}} (E).$  Then $\nabla^{2} \mathcal E=0$ if and only if $(\nabla^{M})^{2} \bar{\mathcal E} =0$, 
 $R^{\nabla^{E}}_{X,\bar{\mathcal E}} =0$, $\mathcal S_{X,  \bar{\mathcal E}} =0$ 
 and
  \begin{equation}\label{cor-E-rel}
 \mathcal S_{X, \nabla^{M}_{Y}\bar{\mathcal E} } s-\frac{1}{2} R^{\nabla^{E}}_{X, \nabla^{M}_{X} \bar{\mathcal E}} s
 +( \Delta_{\mathcal E} - \nabla^{E}_{\bar{\mathcal E}} )( \mathcal S_{X , Y} s -\frac{1}{2} R^{\nabla^{E}}_{X, Y} s) 
 =0
 \end{equation}
 for any $X, Y \in {\mathfrak X}(M)$ and $s\in \Gamma (E)$.
 \end{cor}

 Let $(x^{i})$ be local coordinates on $M$, $( s_{i} )$ a local trivialization of $E$ and $(x^{i}, \xi^{j})$ the induced local coordinates on $E$. 
 Using that  $\nabla$ is torsion-free,  from Lemma 6.2 (6) of \cite{pul} we obtain that 
 \begin{align}
 \nonumber& \nabla_{\partial x^{i}}(\partial x^{j} ) = \nabla^{M} _{\partial x^{i} }( \partial x^{j} ) + \xi^{\alpha} 
 ( \mathcal G_{\partial x^{i},
 \partial x^{j} } s_{\alpha} )^{\uparrow}\\
 \label{pul-nabla}& \nabla_{\partial x^{i} } ( \partial \xi ^{j} ) = \nabla_{\partial \xi^{j}} 
 (\partial x^{i} )
  = (\nabla^{E}_{\partial x^{i}} s_{j}) ^{\uparrow},
  \end{align} 
 and  all other covariant derivatives of coordinate vector fields in direction of coordinate vector fields vanish.
From (\ref{pul-nabla})  we obtain:
\begin{lem}\label{curv-comp} Let $R^{\nabla}$ be the curvature of $\nabla .$ Then
\begin{align}
\nonumber R^{\nabla}_{\partial x^{i}, \partial x^{j} } (\partial x^{k}) &=
R^{\nabla^{M}} _{\partial x^{i}, \partial x^{j} } (\partial x^{k} ) +
u^{\alpha} \left(  \mathcal G_{\partial x^{i}, \nabla^{M}_{\partial x^{j} }\partial x^{k} )}  
s_{\alpha}   +
\nabla^{E}_{\partial x^{i} } (\mathcal G_{\partial x^{j}, \partial x^{k} } s_{\alpha} )\right)\\
\nonumber& - u^{\alpha} \left(  \mathcal G_{\partial x^{j}, \nabla^{M}_{\partial x^{i} }(\partial x^{k} )}  
s_{\alpha}   +
\nabla^{E}_{\partial x^{j} } (\mathcal G_{\partial x^{i}, \partial x^{k}}  s_{\alpha} )\right)\\
\nonumber R^{\nabla}_{\partial x^{i}, \partial u^{j} } (\partial x^{k }) &
=\left( ( \nabla^{E})^{2}_{\partial x^{i}, \partial x^{k} } s_{j} - \mathcal G_{\partial x^{i}, \partial x^{k}} s_{j} \right)^{\uparrow}\\
\nonumber R^{\nabla}_{\partial x^{i}, \partial x^{j} } (\partial u^{k }) &=
\left( R^{\nabla^{E}}_{\partial x^{i}, \partial x^{j}} s_{k} \right)^{\uparrow}
\end{align}
and $R^{{\nabla}}_{X, Y}Z =0$ for all other coordinate vector fields $(X, Y, Z)$. 
\end{lem}
  
\begin{cor}\label{flat-fwl}   The connection $\nabla$ is flat if and only if $\nabla^{M}$, $\nabla^{E}$ are flat and $\mathcal S =0.$
\end{cor}

\begin{proof} From Lemma \ref{curv-comp},  $\nabla$ is flat if and only if $\nabla^{E}$ and $\nabla^{M}$ are flat and 
\begin{align}
\nonumber& \mathcal G_{\partial x^{i}, \nabla^{M}_{\partial x^{j} }(\partial x^{k} )}  
s_{\alpha}   +
\nabla^{E}_{\partial x^{i} } (\mathcal G_{\partial x^{j}, \partial x^{k}}  s_{\alpha} )
-  \mathcal G_{\partial x^{j}, \nabla^{M}_{\partial x^{i} }(\partial x^{k} )}  
s_{\alpha}   - 
\nabla^{E}_{\partial x^{j} } (\mathcal G_{\partial x^{i}, \partial x^{k} } s_{\alpha}) =0,\\
\label{rel-flat}& ( \nabla^{E})^{2}_{\partial x^{i},\partial x^{k} } s_{j}  = 
\mathcal G_{\partial x^{i}, \partial x^{k}} s_{j} ,
\end{align}
for any $ i, j, k$ and $\alpha .$ 
We rewrite the  second relation  in (\ref{rel-flat}) as 
$\mathcal S_{\partial x^{i},\partial x^{j} } s_{\alpha} = \frac{1}{2} R^{\nabla^{E}}_{\partial x^{i} ,
\partial x^{k}} s_{\alpha}$ and obtain  $\mathcal S =0$, since $\nabla^{E}$ is flat.  Using that $\mathcal S =0$, 
$\nabla^{M}$, $\nabla^{E}$ are flat and $\nabla^{M}$ is torsion-free, we obtain the first relation in
(\ref{rel-flat}). 
\end{proof}

\section{Components of FWL $F$-manifolds}\label{components-section}
 
Let  $\pi : E \rightarrow M$ be a vector bundle,   $T\in {\mathcal T}_{\mathrm{fwl},s}^{2,1}(E)$
with components
$(D, l, T^{M})$ 
and  $e\in {\mathfrak X}_{\mathrm{fwl}}(E)$  with  derivation $\Delta_{e}: \Gamma (E) \rightarrow \Gamma (E)$.
We identify $T$  and  $T^{M}$ with 
(commutative) multiplications 
$\circ$ and $*$ on $TE$ and $TM$ respectively. 
 \subsection{FWL multiplication and unit field}\label{components-Fman}

\begin{lem}\label{lem-asoc}
The multiplication   $\circ$ is  associative if and only if
$*$ is associative, 
for any $s\in \Gamma (E)$ and $X, Y\in {\mathfrak X}(M)$, 
\begin{equation}\label{l-xy}
l_{X}  ( l_{Y}s )  = l_{X * Y} s,
\end{equation}
and  the expression
\begin{equation}\label{completely-symm}
l_{Z} (D_{X, Y} s  ) + D_{X * Y, Z} s
\end{equation}
is symmetric in $X, Z\in {\mathfrak X} (M).$ 
\end{lem}

\begin{proof} 
Since $\circ$ is commutative, it is  also associative if and only if 
\begin{equation}\label{asoc-XYZ}
(X\circ Y) \circ Z =  (Z\circ Y)\circ X,
\quad \forall X, Y, Z\in {\mathfrak X}(E).
\end{equation}
Locally, $\circ$ is of the form (\ref{loc}) with $a_{jk}^{i} = \tilde{a}^{i}_{kj}$, $b_{ij}^{a} = b^{a}_{ji}$ and $a^{i}_{j, kp} = a^{i}_{j, pk}.$ 
Relation (\ref{asoc-XYZ}) with two vertical arguments is satisfied (both terms vanish) and 
with one vertical argument is equivalent to
(\ref{l-xy}). 
We now consider (\ref{asoc-XYZ}) with $X= \partial x^{p}$, $Y= \partial x^{q}$ and $Z=\partial x^{r}$.
Projecting (\ref{asoc-XYZ}) to $TM$ we obtain that $*$ is associative. 
The vertical part of the left hand side in 
(\ref{asoc-XYZ})  equals 
\begin{equation}\label{expr-assoc}
\sum_{i, j}\xi^{j} a^{i}_{j, pq}  l_{\partial x^{r}}  s_{i} +\sum_{k, i, j} b^{k}_{pq}  \xi^{j} a^{i}_{j, kr} s_{i}.
\end{equation}
We obtain that (\ref{asoc-XYZ}) holds if and only if 
$$
\sum_{i, j}  \xi^{j} s_{i}^{*} (D_{\partial x^{p}, \partial x^{q}} s_{j}) l_{\partial x^{r}} s_{i} +\sum_{ k, i, j} dx^{k} (\partial x^{p}*  \partial x^{q})
\xi^{j} s_{i}^{*} (D_{\partial x^{k}, \partial x^{r}} s_{j}) s_{i} 
$$
is symmetric in $(r, p)$, or  expression  (\ref{completely-symm}) is symmetric in $(X, Z)$.
\end{proof}

\begin{lem}\label{lem-unit}  Assume that $\circ$ is associative. 
Then  $e$ is a  unit  field for $\circ$ if and only if $\bar{e} = \pi_{*} e$ is a unit field for $*$ and
\begin{equation}\label{e-unit-cond}
l_{\bar{e}} s= s,\ l_{X} (\Delta_{e}s) = D_{\bar{e}  , X} s,
\end{equation}
for any $X\in {\mathfrak X}(M)$ and $s\in \Gamma (E).$
\end{lem}

\begin{proof}
The vector field $e$ is a unit  field  for $\circ$ if and only if 
\begin{equation}\label{e-unit}
e\circ X = X,\ \forall X\in TE.
\end{equation}
Relation (\ref{e-unit}) with
$X:= s^{\uparrow}\in {\mathfrak X}_{\mathrm{core}}(E)$  implies the first relation in  (\ref{e-unit-cond}), as $e\circ s^{\uparrow} =  (l_{\bar{e} } s)^{\uparrow}$. 
As in the previous lemma, $\circ$ is locally of the form (\ref{loc}). Write 
$e= \lambda_{ij} \xi^{i} \partial \xi^{j} +\beta_{i}\partial x^{i}$ where $\lambda_{ij}, \beta_{i}\in C^{\infty}(M)$. 
Consider  (\ref{e-unit}) with
$X:= \partial x^r.$   Projecting it to $M$ we obtain that $\bar{e}$ is a unit  field for $*.$   The  vertical part of $e\circ \partial x^{r}$  is given by 
 \begin{equation}
( \beta_{i} \xi^{j} a^{f}_{j, ir} +\lambda_{ji} \xi^{j} a^{f}_{ir} ) \partial \xi^{f}=
\xi^{j}\left(  \beta_{i} s_{f}^{*} (D_{\partial x^{i}, \partial x^{r}} s_{j}) +\lambda_{ji} s_{f}^{*} (l_{\partial x^{r}} s_{i}) \right)
\partial \xi^{f} 
 \end{equation}
and vanishes if and only if  $D_{\bar{e}, \partial x^{r}} s_{j} = - l_{\partial x^{r}}  (\lambda_{ji }s_{i}).$ 
The 
second relation in (\ref{e-unit-cond}) follows by noticing that 
$\Delta_{e} s_{j} = -\lambda_{ji} s_{i}$.
\end{proof}

\begin{lem}\label{prop-integr} 
Assume that $\circ$ is associative with unit field $e$. Then $(E, \circ , e)$ is an $F$-manifold if and only if 
$(M, * , \bar{e})$ is an $F$-manifold, 
\begin{equation}\label{integr-1}
[ D_{X, Y} , l_{Z}]  s = l_{\mathcal L_{Z} (*)(X, Y)}s 
\end{equation}
and
\begin{align}
\nonumber&   [ D_{Z, V}, D_{X, Y} ]s  \\
\nonumber& = D_{\mathcal L_{X * Y} Z, V} s + D_{\mathcal L_{X * Y} V, Z} s + 
D_{\mathcal L_{X}(*) ( Z, V), Y} s  + D_{\mathcal L_{Y}(*) ( Z, V), X} s \\
\label{F-man-lin} &  - l_{X} ( D_{\mathcal L_{Y}V, Z} s + D_{\mathcal L_{Y}Z, V} s ) -  l_{Y} ( D_{\mathcal L_{X}V, Z} s+  D_{\mathcal L_{X}Z, V} s ),
\end{align}
for any $X, Y, Z, V\in {\mathfrak X}(M)$ and $s\in \Gamma (E)$, where  
$[A, B]= AB - BA$ denotes the commutator of $A, B\in \mathrm{End}\, \Gamma (E).$
\end{lem}

\begin{proof} Define $\mathcal P \in \mathcal T^{4,1}(E)$ by
$$
\mathcal P (X, Y, Z, V) := \mathcal L_{X\circ Y} (\circ )(Z, V)  - X\circ \mathcal L_{Y} (\circ )(Z, V)  - Y\circ \mathcal L_{X} (\circ )(Z, V),
$$
for any $X, Y, Z, V\in {\mathfrak X}(E)$. It satisfies
\begin{align}
\nonumber& \mathcal P (X, Y, Z, V) = \mathcal P (Y, X, Z, V)\\
\nonumber& \mathcal P (X, Y, Z, V) = \mathcal P (X, Y, V, Z)\\
\label{symmetries}& \mathcal P (X, Y, Z, V) = -  \mathcal P (Z, V, X, Y).
\end{align}
From Lemma \ref{lin-core}, 
$\mathcal P$ is FWL.  Its basic component is given by
\begin{equation}\label{r}
\tilde{r}(\bar{X}, \bar{Y}, \bar{Z}, \bar{V}) = \mathcal L_{\bar{X} * \bar{Y}} ( * )(\bar{Z}, \bar{V})  - 
\bar{X} * \mathcal L_{\bar{Y}} ( * )(\bar{Z}, \bar{V})  - \bar{Y}* \mathcal L_{\bar{X}} (*)(\bar{Z}, \bar{V}).
\end{equation}
Let  $(\tilde{D}, \tilde{l}^{(i)})$  
be the remaining components of $\mathcal P$,  where 
\begin{align}
\nonumber& \tilde{D} : \Gamma (E) \rightarrow \Gamma (  (T^{*}M)^{4} \otimes E)\\
\nonumber&\tilde{ l}^{(i)} : E \rightarrow ( T^{*}M)^{3} \otimes E,\ 1\leq i\leq 4.
\end{align}
We claim that 
\begin{align}
\nonumber& \tilde{l}^{(1)} = \tilde{l}^{(2)},\ \tilde{l}^{(3)} = \tilde{l}^{(4)},\   \tilde{l}^{(3)}(s)( \bar{X}, \bar{Y}, \bar{Z}) = - \tilde{l}^{(1)} (s)(\bar{Z}, \bar{X}, \bar{Y})\\
\label{com-l}& \tilde{l}^{(1)} (s)( \bar{X}, \bar{Y}, \bar{Z}) = [ D_{\bar{Y}, \bar{Z}} , l_{\bar{X}} ] s - l_{\mathcal L_{\bar{X}} (* ) (\bar{Y}, \bar{Z})}s,
\end{align}
for any $\bar{X}, \bar{Y}, \bar{Z}, \bar{V}\in  {\mathfrak X}(M)$ and $s\in \Gamma (E).$
The first line in (\ref{com-l}) follows from  (\ref{symmetries}). 
 For the second line  in  (\ref{com-l}), let  $Y, Z, V\in  {\mathfrak X}(E)$ projectable. Then
 \begin{align}
\nonumber&  \mathcal L_{s^{\uparrow}\circ Y} (\circ ) (Z, V) = \mathcal L_{(l_{\bar{Y}} s)^{\uparrow}}(\circ )(Z, V) = D_{\bar{Z}, \bar{V}}( l_{\bar{Y}} s )^{\uparrow}\\
\nonumber& s^{\uparrow} \circ\mathcal L_{Y} (\circ) (Z, V) =\left(  l_{\mathcal L_{\bar{Y}} (*) (\bar{Z}, \bar{V})} s\right)^{\uparrow}\\
\nonumber&Y\circ {\mathcal L}_{s^{\uparrow}}  (\circ )  (Z, V) = Y\circ (D_{\bar{Z}, \bar{V}} s)^{\uparrow}
= l_{\bar{Y}} (D_{\bar{Z}, \bar{V}}s)^{\uparrow} 
\end{align}
which imply 
\begin{equation}
\tilde{l}^{(1)} (s) (\bar{Y}, \bar{Z}, \bar{V})^{\uparrow} 
= \mathcal P (s^{\uparrow}, Y, Z ,V) = \left( D_{\bar{Z}, \bar{V}}( l_{\bar{Y}} s )- l_{\bar{Y}} (D_{\bar{Z}, \bar{V}}s)
 - l_{\mathcal L_{\bar{Y}} (*) (\bar{Z}, \bar{V})} s\right)^{\uparrow}
\end{equation}
i.e. the second line in (\ref{com-l}) holds.  We now prove  that 
\begin{align}
(\tilde{D}s)( \bar{X}, \bar{Y}, \bar{Z}, \bar{V}) \nonumber& = - D_{\mathcal L_{\bar{X} * \bar{Y}}\bar{Z}, 
\bar{V}} s - D_{\mathcal L_{\bar{X} * \bar{Y}}\bar{V}, \bar{Z}} s  - D_{\mathcal L_{\bar{X}}(* )( \bar{Z}, \bar{V}), \bar{Y}} s\\
\nonumber& - D_{\mathcal L_{\bar{Y}}( * )( \bar{Z}, \bar{V})} s
 + l_{\bar{X}} ( D_{\mathcal L_{\bar{Y}}\bar{V}, \bar{Z}} s + D_{\mathcal L_{\bar{Y}}\bar{Z}, \bar{V}} s ) \\
\label{tilde-d}&  + l_{\bar{Y}} ( D_{\mathcal L_{\bar{X}}\bar{V},\bar{ Z}} s + D_{\mathcal L_{\bar{X}}\bar{Z}, \bar{V}} s ) - [ D_{\bar{X}, \bar{Y}}, 
D_{\bar{Z}, \bar{V}} ]s,
\end{align}
for any $\bar{X}, \bar{Y}, \bar{Z}, \bar{V}\in \mathfrak {X} (M).$ For any  $X, Y, Z, V \in {\mathfrak X}_{\mathrm{fwl}} (E)$, 
\begin{align}
&(\nonumber {\mathcal L_{s^{\uparrow}} \mathcal P} )(X, Y, Z, V)  = \mathcal L_{s^{\uparrow}} ( \mathcal P (X, Y, Z, V)) - \mathcal P (\mathcal L_{s^{\uparrow}}
X, Y, Z, V) \\
\label{p}& -  \mathcal P (X, \mathcal L_{s^{\uparrow}} Y, Z, V)  -  \mathcal P ( X, Y, \mathcal L_{s^{\uparrow}}Z, V) -
 \mathcal P ( X, Y,  Z, \mathcal L_{s^{\uparrow}}V) .
\end{align}
Expanding the first term in the right hand side of (\ref{p}), 
\begin{align}
\nonumber& \mathcal L_{s^{\uparrow}} ( \mathcal P (X, Y, Z, V)) = \mathcal L_{s^{\uparrow}} (\mathcal L_{X\circ Y} (Z\circ V))
- \mathcal L_{s^{\uparrow}} (\circ ) (\mathcal L_{X\circ Y}Z, V) - \mathcal L_{s^{\uparrow}} ( \mathcal L_{X\circ Y}Z) \circ V\\
\nonumber& - (\mathcal L_{X\circ Y} Z)\circ \mathcal L_{s^{\uparrow}} V - \mathcal L_{s^{\uparrow}} (\circ ) (Z, \mathcal L_{X\circ Y}V)
-(\mathcal L_{s^{\uparrow}} Z)\circ \mathcal L_{X\circ Y} V - Z\circ \mathcal L_{s^{\uparrow}} \mathcal L_{X\circ Y} V\\
\nonumber& - \mathcal L_{s^{\uparrow}} (\circ ) (X, \mathcal L_{Y}(Z\circ V) -( \mathcal L_{Y} Z)\circ V - Z\circ \mathcal L_{Y} V)\\
\nonumber& - (\mathcal L_{s^{\uparrow}} X)\circ \left( \mathcal L_{Y}(Z\circ V) -( \mathcal L_{Y} Z)\circ V - Z\circ \mathcal L_{Y} V)\right)\\
\nonumber& - X\circ \mathcal L_{s^{\uparrow}}  \left( \mathcal L_{Y}(Z\circ V) -( \mathcal L_{Y} Z)\circ V - Z\circ \mathcal L_{Y} V\right)\\
\nonumber& - \mathcal L_{s^{\uparrow}} (\circ ) (Y, \mathcal L_{X}(Z\circ V) -( \mathcal L_{X} Z)\circ V - Z\circ \mathcal L_{X} V)\\
\nonumber& - (\mathcal L_{s^{\uparrow}} Y)\circ \left( \mathcal L_{X}(Z\circ V) -( \mathcal L_{X} Z)\circ V - Z\circ \mathcal L_{X} V)\right)\\
\nonumber& - Y\circ \mathcal L_{s^{\uparrow}}  \left( \mathcal L_{X}(Z\circ V) -( \mathcal L_{X} Z)\circ V - Z\circ \mathcal L_{X} V\right) .
\end{align}
Using  repeatedly 
(\ref{deriv-vector-field}), 
 (\ref{data-1}) and  (\ref{data-2}), that  the Lie bracket  of  two FWL vector fields is  FWL
 and $X\circ Y$ is FWL whenever $X$ and $Y$ are FWL,  we obtain 
\begin{align}
\nonumber& \mathcal L_{s^{\uparrow}} ( \mathcal P (X, Y, Z, V))  = - (\Delta_{\mathcal L_{X\circ Y} (Z\circ V)} s)^{\uparrow}
- (D_{\mathcal L_{\bar{X}*\bar{Y}}\bar{Z}, \bar{V}} s)^{\uparrow} + l_{\bar{V}} (\Delta_{\mathcal L_{X\circ Y}Z}s)^{\uparrow}\\
\nonumber&  + l_{\mathcal L_{\bar{X}*\bar{Y}} \bar{Z}} 
(\Delta_{V} s)^{\uparrow} - (D_{\bar{Z}, \mathcal L_{\bar{X}*\bar{Y}}\bar{V}}s)^{\uparrow} 
+( l_{\mathcal L_{\bar{X}*\bar{Y}} \bar{V}} \Delta_{Z} s)^{\uparrow}  + l_{\bar{Z}} (\Delta_{\mathcal L_{X\circ Y} V}s)^{\uparrow} \\
\nonumber& - (D_{\bar{X}, \mathcal L_{\bar{Y}} (*) (\bar{Z}, \bar{V}}s)^{\uparrow} + l_{\mathcal L_{\bar{Y}}(*) (\bar{Z}, \bar{V})} (\Delta_{X} s)^{\uparrow} +( l_{\bar{X}} ( \Delta_{\mathcal L_{Y}(\circ ) (Z, V)} s)^{\uparrow}\\
\nonumber& - (D_{\bar{Y}, \mathcal L_{\bar{X}} (*) (\bar{Z}, \bar{V}}s)^{\uparrow} + l_{\mathcal L_{\bar{X}}(*) (\bar{Z}, \bar{V})} (\Delta_{Y} s)^{\uparrow} +( l_{\bar{Y}} ( \Delta_{\mathcal L_{X}(\circ ) (Z, V)} s)^{\uparrow} 
\end{align}
By similar computations  we obtain the other  terms from the right hand side of (\ref{p}).  Using (\ref{comutator}) and cancelling terms 
we  obtain that $(\nonumber {\mathcal L_{s^{\uparrow}} \mathcal P} )(X, Y, Z, V)$ is equal to
the core vector field defined by the right hand side of 
(\ref{tilde-d}). Relation (\ref{tilde-d}) follows.  Relations (\ref{r}), (\ref{com-l})  and (\ref{tilde-d}) imply our claim.
\end{proof}

\begin{defn}\label{over} Let  $\pi : E \rightarrow M$ be a vector bundle and $(E, \circ , e )$ a FWL $F$-manifold.
Let $*\in {\mathcal T}^{2,1} (M)$ be the basic component of $\circ $ and $\bar{e} := \pi_{*} e$.
Then $(E, \circ , e)$ is called a  {\cmssl  FWL $F$-manifold over  the $F$-manifold  $(M, *,\bar{e})$.}
\end{defn}

\subsection{FWL Euler fields}\label{euler-sect}

Let $(M, \circ , e)$ be a FWL $F$-manifold and $(D, l, *)$ the components of $\circ .$ 

\begin{prop}\label{euler-cor} A FWL  vector field $\mathcal E\in {\mathfrak X}_{\mathrm{fwl}} (E)$ is Euler on
$(E, \circ , e)$ if and only if 
$\bar{\mathcal E} =  \pi_{*} \mathcal E$ is Euler on $(M, *, \bar{e})$ and,  for any $X, Y\in {\mathfrak X}(M)$,
\begin{align}
\nonumber& [ \Delta_{\mathcal E}, l_{X}]  -  l_{\mathcal L_{\bar{\mathcal E}} X}= l_{X} \\
\label{euler-cond}& [ \Delta_{\mathcal E}, D_{X, Y} ]  -  D_{\mathcal L_{\bar{\mathcal E}}X, Y}   -  D_{Y, \mathcal L_{\bar{\mathcal E}}Y} 
= D_{X, Y} .
\end{align}
\end{prop}

\begin{proof}  
Compare   the components of $\mathcal L_{\mathcal E}(\circ )$ and  $\circ$ using   Lemma \ref{prelim} below. 
\end{proof}

\begin{lem}\label{prelim} Let  $\mathcal E\in {\mathfrak X}_{\mathrm{fwl}} (E)$. Then $\mathcal L_{\mathcal E } (\circ )\in
\mathcal T^{2,1}_{\mathrm{fwl},s} (E)$ has  components
$(\tilde{D}, \tilde{l}, \tilde{r})$ where 
$\tilde{r} = \mathcal L_{\bar{\mathcal E}} (*)$, and, for any $X, Y\in {\mathfrak X}(M)$, 
\begin{align}
\nonumber& \tilde{D}_{X, Y}  =  [ \Delta_{\mathcal E}  , D_{X, Y} ]  - D_{\mathcal L_{\bar{\mathcal E}}X, Y}  - D_{X, \mathcal L_{\bar{\mathcal E}}Y} \\
\label{comp-lie}& \tilde{l}_{X}  =[ \Delta_{\mathcal E},  l_{X} ]s - l_{\mathcal L_{\bar{\mathcal E}} X} .
\end{align}
\end{lem}

\begin{proof} 
From Lemma \ref{lin-core},  $\mathcal L_{\mathcal E} (\circ )$ is FWL. 
For any $s\in \Gamma (E)$,  $X\in {\mathfrak X}_{\mathrm{fwl}} (E)$, 
\begin{align}
\nonumber&  \mathcal L_{\mathcal E} ( \circ ) (s^{\uparrow}, X) = \mathcal L_{\mathcal E} ( s^{\uparrow} \circ X) -   (\mathcal L_{\mathcal E} s^{\uparrow}) \circ X - s^{\uparrow} \circ \mathcal L_{\mathcal E} X \\
\nonumber& =\mathcal L_{\mathcal E} (  (l_{\bar{X}} s)^{\uparrow} ) - ( \Delta_{\mathcal E} s)^{\uparrow} \circ X - ( l_{\mathcal L_{\bar{\mathcal E}} \bar{X}}s)^{\uparrow}= \left( \Delta_{\mathcal E} ( l_{\bar{X}} s) - l_{\bar{X}} (\Delta_{\mathcal E} s)-  l_{\mathcal L_{\bar{\mathcal E}} \bar{X} } s\right)^{\uparrow}
\end{align}
which implies the second  relation  in
(\ref{comp-lie}).  The identity
$$
 \mathcal L_{s^{\uparrow}} \mathcal L_{\mathcal E} (\circ )  
= \mathcal L_{\mathcal E} \mathcal L_{s^{\uparrow}} (\circ ) + \mathcal L_{\mathcal L_{s^{\uparrow}} \mathcal E } (\circ )
$$
implies,  for any $X, Y\in {\mathfrak X}_{\mathrm{fwl}} (E)$, 
$$
(\tilde{D}_{\bar{X}, \bar{Y}} s)^{\uparrow} = ( \mathcal L_{s^{\uparrow}}  \mathcal L_{\mathcal E} (\circ ))  (X, Y)= (\mathcal L_{\mathcal E}  \mathcal L_{s^{\uparrow}} (\circ ))  (X, Y) -\mathcal L_{ (\Delta_{\mathcal E} s)^{\uparrow} } (\circ ) (X, Y).
$$
But $\mathcal L_{ (\Delta_{\mathcal E} s)^{\uparrow} } (\circ ) (X, Y) = D_{\bar{X}, \bar{Y}} (\Delta_{\mathcal E} s)^{\uparrow} $ 
and 
\begin{align}
\nonumber&  ( \mathcal L_{\mathcal E} \mathcal L_{s^{\uparrow}} (\circ ) ) (X, Y) = \mathcal L_{\mathcal E} \left( \mathcal L_{s^{\uparrow}} (\circ ) (X, Y) \right) -
\mathcal L_{s^{\uparrow}} (\circ ) (\mathcal L_{\mathcal E} X, Y) - \mathcal L_{s^{\uparrow}} (\circ ) (X, \mathcal L_{\mathcal E}  Y) \\
\nonumber& = \left( \Delta_{\mathcal E} ( D_{\bar{X}, \bar{Y}} s) - D_{\mathcal L_{\bar{\mathcal E}} \bar{X}, \bar{Y} } s
- D_{\bar{X}, \mathcal L_{\bar{\mathcal E}} \bar{Y} }s\right)^{\uparrow } .
\end{align}
The first relation in (\ref{comp-lie}) follows. Relation $\tilde{r} =\mathcal L_{\bar{\mathcal E}} (*)$ is obvious.
\end{proof}
\begin{defn}\label{over-euler} Let  $\pi : E \rightarrow M$ be a vector bundle and $(E, \circ , e , \mathcal E)$ a  FWL $F$-manifold
with Euler field.
Let $*\in {\mathcal T}^{2,1} (M)$ be the basic component of $\circ $ and $\bar{e} := \pi_{*} e$,   $\bar{\mathcal E} = \pi_{*} \mathcal E $.
Then $(E, \circ , e,\mathcal E )$ is called a  {\cmssl  FWL $F$-manifold with Euler field over  the $F$-manifold  with Euler field $(M, *,\bar{e},
\bar{\mathcal E})$.}
\end{defn}

\subsection{FWL compatible connections}\label{conn-sect}

Let $(E, \circ , e)$ be a  FWL $F$-manifold  
and $\nabla$ a FWL torsion-free connection on $E$. 
Let  $(D, l, *)$  and $(\nabla^{M}, \nabla^{E}, \mathcal S )$ be the components of  $\circ $ and   $\nabla$ respectively
and   $\mathcal G$  defined by 
(\ref{g-s-1}).

\begin{prop}\label{nabla-circ-sym} The tensor-field $\nabla (\circ )$ is totally symmetric if and only if $\nabla^{M} (*)$ is totally symmetric and,  for any
$X, Y, Z\in {\mathfrak X}(M)$ and $s\in \Gamma (E)$,
\begin{align}
\nonumber& \mathcal G_{X, Y*Z} s - \mathcal G_{Y, X*Z} s+\nabla^{E}_{X} (D_{Y, Z}s) -\nabla^{E}_{Y} (D_{X, Z} s)\\
\label{con-1}&  - l_{Y} \mathcal G_{X, Z}s + l_{X} \mathcal G_{Y, Z} s - D_{Y, \nabla^{M}_{X} Z}s + D_{X, \nabla^{M}_{Y} Z} s 
- D_{\mathcal L_{X}Y, Z} s =0
\end{align}
and
\begin{equation}\label{con-2}
\nabla^{E}_{Y*Z} s+ D_{Y, Z} s - l_{Y} \nabla^{E}_{Z} s -\nabla^{E}_{Y} (l_{Z} s) + l_{\nabla^{M}_{Y} Z} s =0
\end{equation}
and
\begin{equation}\label{con-3}
\nabla^{E}_{X}( l_{Y} s) -\nabla^{E}_{Y} ( l_{X} s) - l_{\mathcal L_{X} Y}s + l_{X} \nabla^{E}_{Y} s - l_{Y} \nabla^{E}_{X} s =0.
\end{equation}
\end{prop}

\begin{proof} Define
$\mathcal P\in {\mathcal T}^{3,1}(E)$ by 
$\mathcal P (X, Y, Z) := \nabla_{X}(\circ )(Y, Z) - \nabla_{Y} (\circ )(X, Z)$. It  is FWL  since $\mathcal P (X, Y, Z)$ is FWL whenever
$X, Y, Z$ are so. Its basic component equals 
$(\bar{X}, \bar{Y}, \bar{Z}) \mapsto \nabla^{M}_{\bar{X}}(*)(\bar{Y}, \bar{Z}) -\nabla^{M}_{\bar{Y}}(*)(\bar{X}, \bar{Z}).$ 
Let $X, Y, Z\in {\mathfrak X}_{\mathrm{\mathrm{fwl}}} (E).$  Then
\begin{align*}
\nonumber& \mathcal P (s^{\uparrow}, Y, Z)=  \nabla_{s^{\uparrow}}(\circ )(Y, Z) - \nabla_{Y} (\circ )(s^{\uparrow}, Z)=\\
\nonumber&\nabla_{s^{\uparrow}} (Y\circ Z) -(\nabla_{s^{\uparrow}} Y)\circ Z - Y\circ (\nabla_{s^{\uparrow}} Z ) 
-\nabla_{Y} (s^{\uparrow}\circ Z) +(\nabla_{Y} s^{\uparrow}) \circ Z + s^{\uparrow}\circ \nabla_{Y} Z\\
\nonumber&  =(\nabla^{E}_{\bar{Y}*\bar{Z}} s -\Delta_{Y\circ Z} s
+ l_{\bar{Z}} \Delta_{Y} s +l_{\bar{Y}}( \Delta_{Z} s - \nabla^{E}_{\bar{Z}})  s -\nabla^{E}_{\bar{Y}} ( l_{\bar{Z}} s) + l_{\nabla^{M}_{\bar{Y}}\bar{Z}} s)^{\uparrow},
\end{align*}
where we used (\ref{fwl-E}) and the first relation in (\ref{need-nabla}).
We claim that 
\begin{equation}\label{nabla-xy-D}
\Delta_{Y\circ Z} s = - D_{\bar{Y}, \bar{Z}} s +  l_{\bar{Y}} \Delta_{Z} s +  l_{\bar{Z}} \Delta_{Y} s.
\end{equation}
This follows from the next computation:
\begin{align}
\nonumber& (\Delta_{Y\circ Z} s)^{\uparrow} = \mathcal L_{Y\circ Z} (s^{\uparrow}) = - \mathcal L_{s^{\uparrow}} (\circ )(Y, Z) -( {\mathcal L}_{s^{\uparrow}} Y)\circ Z - Y\circ \mathcal L_{s^{\uparrow}}Z\\
\nonumber&  =  - (D_{\bar{Y}, \bar{Z}} s)^{\uparrow} +( \Delta_{Y} s)^{\uparrow} \circ Z + Y \circ (\Delta_{Z} s)^{\uparrow}  =  ( - D_{\bar{Y}, \bar{Z}} s +l_{\bar{Z}} ( \Delta_{Y} s) + l_{\bar{Y}} \Delta_{Z} s)^{\uparrow}.
\end{align}
Replacing  $\Delta_{Y\circ Z}s$  with the right hand side of (\ref{nabla-xy-D}) 
we obtain that   $\mathcal P (s^{\uparrow}, Y, Z)$ is the core vector field defined by the section  of $E$ from the left hand side of 
(\ref{con-2}) with $Y, Z$ replaced by $\bar{Y}, \bar{Z}$. 
Similarly,  $\mathcal P (X, Y, s^{\uparrow}) $ is the core vector field defined by the section  of $E$ from the left hand side of (\ref{con-3}) with $X, Y$ replaced by $\bar{X}, \bar{Y}$. 
Using in addition 
(\ref{deriv-G}) and  the second relation in  (\ref{need-nabla}), we obtain that $(\mathcal L_{s^{\uparrow}} \mathcal P  )(X, Y, Z)$
is the core vector field defined by the section  of $E$ from  the left hand side of (\ref{con-1}) with $X, Y, Z$ replaced by $\bar{X}, \bar{Y}, \bar{Z}$.  Using $\mathcal P ( s^{\uparrow}, Y, Z) = - \mathcal P (Y, s^{\uparrow}, Z)$ and that a FWL tensor field vanishes if and only if all its components vanish, we obtain our claim. 
\end{proof}

\begin{lem}\label{mathcal-s} Assume  (\ref{con-2}) and (\ref{con-3}) hold.  Then  (\ref{con-1}) is equivalent to
\begin{align}
\nonumber&\mathcal S_{X, Y*Z} s- \mathcal S_{Y, X*Z} s+ l_{X} \mathcal S_{Y, Z}s - l_{Y} \mathcal S_{X, Z} s + R^{\nabla^{E}}_{X, Y} (l_{Z} s) - l_{R^{\nabla^{M}}_{X, Y} Z}s\\
\label{con-1-alt}& = \frac{1}{2} \left( R^{\nabla^{E}}_{X, Y*Z}s - R^{\nabla^{E}}_{Y, X*Z}s  + l_{X} R^{\nabla^{E}}_{Y, Z}s- l_{Y} R^{\nabla^{E}}_{X, Z}s \right)  
\end{align}
for any $X, Y, Z\in {\mathfrak X}(M)$  and $s \in \Gamma (E).$ 
\end{lem}

\begin{proof} In terms of  the $\mathcal S$-component   of $\nabla$, relation   (\ref{con-1}) becomes
\begin{align}
\nonumber& \mathcal S_{X, Y*Z} s - \mathcal S_{Y, X* Z} s 
+ l_{X} \mathcal S_{Y, Z} s -  l_{Y} \mathcal S_{X, Z} s \\
\nonumber& - \frac{1}{2} \left(  R^{\nabla^{E}}_{X , Y* Z} s - R^{\nabla^{E}}_{Y , X*Z} s 
+ l_{X} R^{\nabla^{E}} _{Y, Z} s   - l_{Y} R^{\nabla^{E}} _{X, Z} s \right)\\
\nonumber & + (\nabla^{E})^{2}_{X, Y*Z} s -(\nabla^{E})^{2}_{Y, X*Z} s + l_{X} (\nabla^{E})^{2}_{Y, Z}s- l_{Y} (\nabla^{E})^{2}_{X, Z} s \\
\label{con-alt} & +\nabla^{E}_{X} (D_{Y, Z}s) -\nabla^{E}_{Y} (D_{X, Z} s)  + D_{X, \nabla^{M}_{Y}Z}s - D_{Y,\nabla^{M}_{X}Z} s 
- D_{\mathcal L_{X}Y, Z} s=0.
\end{align}
From (\ref{con-2}), the  sum of the first terms in the third and fourth lines equals 
$$
 (\nabla^{E})^{2}_{X, Y*Z} s +\nabla^{E}_{X} (D_{Y, Z}s)=\nabla^{E}_{X} \left(l_{Y} \nabla^{E}_{Z} s +\nabla^{E}_{Y} (l_{Z} s) -l_{\nabla^{M}_{Y}Z} s
 \right) -\nabla^{E}_{\nabla^{M}_{X} (Y*Z)} s.
 $$
Computing similarly $(\nabla^{E})^{2}_{Y, X*Z} s +\nabla^{E}_{Y} (D_{X, Z}s)$ we obtain that the last two lines in the 
left hand side of (\ref{con-alt}) equal 
\begin{align}
\nonumber&\nabla^{E}_{X} (l_{Y}) (\nabla^{E}_{Z}s) +\nabla^{E}_{X}\nabla^{E}_{Y} (l_{Z}s) -\nabla^{E}_{X}(l_{\nabla^{M}_{Y}Z} s) 
-\nabla^{E}_{\nabla^{M}_{X}(Y*Z)}s - \nabla^{E}_{Y} (l_{X})\nabla^{E}_{Z}s\\
\nonumber& - \nabla^{E}_{Y} \nabla^{E}_{X} (l_{Z}s)+\nabla^{E}_{Y} (l_{\nabla^{M}_{X}Z} s) +\nabla^{E}_{\nabla^{M}_{Y} (X*Z)} s + l_{Y} (\nabla^{E}_{\nabla^{M}_{X}Z} s)
- l_{X} (\nabla^{M}_{\nabla^{M}_{Y} Z}s)\\
\label{S-nonflat}&  - D_{Y,\nabla^{M}_{X}Z}s + D_{X, \nabla^{M}_{Y}Z} s
- D_{\mathcal L_{X}Y, Z}s,
\end{align}
From (\ref{con-3}), 
$\nabla^{E}_{X} (l_{Y}) \nabla^{E}_{Z} s-\nabla^{E}_{Y} (l_{X})\nabla^{E}_{Z} s= l_{\mathcal L_{X}Y} ( \nabla^{E}_{Z}s)$. 
Since $\nabla^{M}$ is compatible,
\begin{equation}
-\nabla^{E}_{\nabla^{M}_{X}(Y*Z)}s+\nabla^{E}_{\nabla^{M}_{Y} (X*Z)} s = \nabla^{E}_{(\mathcal L_{Y} X)* Z  + X*\nabla^{M}_{Y} Z - Y* \nabla^{M}_{X}Z}s.
\end{equation}
Replacing these expressions in (\ref{S-nonflat})  and applying    (\ref{con-2})
for $D_{X, \nabla^{M}_{Y}Z} s +\nabla^{E}_{X*\nabla^{M}_{Y}Z}s$, $D_{Y, \nabla^{M}_{X}Z} s +\nabla^{E}_{Y*\nabla^{M}_{X}Z}s$
and $D_{\mathcal L_{X}Y, Z} s +\nabla^{E}_{(\mathcal L_{X}Y)*Z}s$,
we obtain that the last two lines in the left hand side of
(\ref{con-alt}) equal $R^{\nabla^{E}}_{X, Y} (l_{Z}s) - l_{R^{\nabla^{M}}_{X, Y}Z}s.$
\end{proof}

\begin{prop}\label{mathcal-e} The unit field $e$ is $\nabla$-parallel if and only if $\bar{e}= \pi_{*} e$ is $\nabla^{M}$-parallel and for any 
$s\in \Gamma (E)$ and $X\in {\mathfrak X}(M)$, 
\begin{equation}\label{nabla-fwl-e}
\nabla^{E}_{\bar{e}} s = \Delta_{e} s,\ \mathcal S_{X,\bar{e}} s= \frac{1}{2} R^{\nabla^{E}}_{X, \bar{e}} s.
\end{equation}
\end{prop} 

\begin{proof}  The basic component of  $\nabla e\in \mathcal T^{1,1}_{\mathrm{fwl}}(E)$ is
$\nabla^{M}\bar{e}$. For any $s\in \Gamma (E)$, 
\begin{align}
\nonumber& \nabla_{s^{\uparrow}} e = \nabla_{e} (s^{\uparrow}) +\mathcal L_{s^{\uparrow}} e= (\nabla^{E}_{\bar{e}} s - \Delta_{e} s)^{\uparrow}\\
\label{nabla-fwl-ep} &\mathcal L_{s^{\uparrow}}  (\nabla e)_{X}= \left( \mathcal G_{\bar{X},\bar{e}} s- \nabla^{E}_{\bar{X}} (\Delta_{e} s)  \right)^{\uparrow}.
\end{align}
Therefore, $\nabla e=0$ if and only if $\nabla^{M}\bar{e}=0$, 
the first relation in (\ref{nabla-fwl-e}) holds and 
$\mathcal G_{\bar{X},\bar{e}}s = \nabla^{E}_{\bar{X}} (\Delta_{e} s).$ 
Replacing $\Delta_{e} s$ with $\nabla^{E}_{\bar{e}} s$
and using 
(\ref{g-s-1}) and  $\nabla^{M}\bar{e} =0$ 
we obtain the second relation in (\ref{nabla-fwl-e}) from the second relation in 
(\ref{nabla-fwl-ep}). \end{proof} 

\begin{defn}\label{over-conn}
Let  $\pi : E \rightarrow M$ be a vector bundle and $(E, \circ , e , \nabla)$ a  FWL $F$-manifold
with compatible connection.
Let $*\in {\mathcal T}^{2,1} (M)$  and $\nabla^{M}$ be the basic components of $\circ $ and $\nabla$  respectively and $\bar{e} := \pi_{*} e$.
Then $(E, \circ , e,\nabla  )$ is called a  {\cmssl  FWL $F$-manifold with compatible connection  over  the $F$-manifold  with 
compatible connection  $(M, *,\bar{e},\nabla^{M})$.}
\end{defn}

\section{A duality for  FWL $F$-manifolds}\label{duality-section}

Let  $\pi : E \rightarrow M$ be a vector bundle, 
$(E, \circ , e)$ a FWL $F$-manifold and $\nabla$ a connection on $M$ with torsion $T^{\nabla}.$ 
Let  $(D, l, r = * )$ be the components of $\circ $.
Define  $D^{*} : \Gamma (E^{*}) \rightarrow \Gamma ( S^{2} T^{*}M\otimes E^*)$, 
$l^{*} : E^{*} \rightarrow T^{*}M\otimes  E^{*}$ and $r^{*} : TM\otimes TM\rightarrow TM$ 
by $r^{*} (X, Y)= X*Y$,  $(l^{*} _{X} \mu ) s := \mu ( l_{X} s)$ and 
\begin{equation}\label{d*}
(D^{*}_{X, Y} \mu ) s:= X \left( \mu ( l_{Y} s)\right) + Y \left( \mu ( l_{X} s)\right)  - \mu \left( l_{[X,Y]_{\nabla} } s \right)  - (X*Y) ( \mu s) -\mu ( D_{X, Y} s)
\end{equation}
for any $X, Y\in {\mathfrak X}(M)$, $\mu \in \Gamma (E^{*})$,  $s\in \Gamma (E)$.
From Proposition \ref{comp-prop}, $(D^{*}, l^{*} , r^{*})$ are the components of a  FWL commutative multiplication  $\bar{\circ}\in \mathcal T^{2,1}_{\mathrm{fwl},s} (E^{*}) $.

\begin{lem}\label{duality-asoc} The multiplication $\bar{\circ}$ is associative if and only if
\begin{align}
\nonumber& l(s)\left(  T^{\nabla}_{X, Y} * Z + 2 T^{\nabla}_{X, Z} * Y + T^{\nabla}_{Y, Z} * X \right)\\
\label{asoc}& +  l(s)\left(  T^{\nabla} _{Z, X* Y} - T^{\nabla}_{X, Y* Z} + 2 \nabla_{X}  (*)  (Y, Z) - 2 \nabla_{Z}  (*)(X, Y)  \right) =0,
\end{align}
for any $s\in \Gamma (E)$ and $X, Y, Z\in {\mathfrak X}(M).$
\end{lem}

\begin{proof} We apply  Lemma \ref{lem-asoc}. 
Relation (\ref{l-xy}) for $(D^{*}, l^{*}, r^{*})$ follows from relation (\ref{l-xy})  for
$(D, l, r).$ For the symmetry  condition, let
\begin{equation}
\mathcal E (X, Y, Z, \mu ):= l^{*}_{Z}  (D^{*}_{X, Y} \mu )+ D^{*}_{X* Y, Z} \mu \in \Gamma (E^{*})
\end{equation}
where  $X, Y, Z\in {\mathfrak X} (M)$ and $\mu \in \Gamma (E^{*}).$  
From the definition of $D^{*}$,  $l^{*}$ and 
(\ref{l-xy})  we obtain that  $\mathcal E (X, Y, Z, \mu ) s$ is symmetric in  $(X, Z)$ 
if and only if
\begin{align}
D_{X, Y} (l_{Z} s)+D_{X*Y, Z} s +  l_{[X*Y , Z]_{\nabla} }s    + l_{ Z*[X,Y]_{\nabla}} s, 
\end{align}
or  
\begin{equation}\label{intermediate}
l(s) ( \mathcal L_{Z} (*) (X, Y) + [X*Y,  Z]_{\nabla} +  Z*  [X,Y]_{\nabla}    ) ,
\end{equation}
is symmetric in $(X, Z)$ for any $s\in \Gamma (E)$, where  in (\ref{intermediate}) 
we used that 
$D_{X*Y, Z} s + l_{Z} ( D_{X, Y} s)$
is symmetric in $(X, Z)$ and  
$[D_{X, Y},  l_{Z} ]s  =l_{\mathcal L_{Z} (*)(X, Y)}s$
(see Lemmas \ref{lem-asoc} and \ref{prop-integr}). 
Expanding  in  (\ref{intermediate})  
$$
[ X*Y, Z]_{\nabla} = \nabla_{X*Y}Z +\nabla_{Z} (*) (X, Y) +(\nabla_{Z} X)*Y + X* \nabla_{Z}Y
$$
and 
$$
\mathcal L_{Z} (*) (X, Y)=\mathcal L_{Z} (X*Y) - (\mathcal L_{Z} X)*Y -X* (\mathcal L_{Z}Y)
$$
and using 
\begin{align}
\nonumber& \mathcal L_{Z} (X*Y) =\nabla_{Z} (X*Y) -\nabla_{X*Y} Z - T^{\nabla} (Z, X*Y)\\
\nonumber& = \nabla_{Z} (*) (X, Y) + (\nabla_{Z} X) * Y + X *  (\nabla_{Z} Y)-\nabla_{X*Y} Z - T^{\nabla} (Z, X*Y)
\end{align}
we obtain that  (\ref{intermediate})  is symmetric in $(X, Z)$ if and only if 
(\ref{asoc}) holds.
\end{proof}

\begin{lem}\label{duality-unit}  Assume that $\bar{\circ}$ is associative. Then $e^{*}$ is a unit  field for $\bar{\circ}$ if and only if 
\begin{equation}\label{unit-duality}
l(s)(   2\nabla_{X} \bar{e} + T^{\nabla }_{\bar{e}, X} ) =0,
\end{equation}
for any   $s\in \Gamma (E)$ and $X\in \mathfrak{X}(M).$
\end{lem}

\begin{proof}  
We apply Lemma \ref{lem-unit}.
We only  check that 
\begin{equation}\label{e-star}
l^{*}_{X} (\Delta_{e^{*}} \mu ) = D^{*}_{\bar{e}, X}\mu ,\ \forall \mu\in \Gamma (E^{*}),\ X\in {\mathfrak X}(E)
\end{equation}
as the other conditions  are trivial. Using $l_{\bar{e}} s= s$,  we obtain that (\ref{e-star}) is equivalent to 
\begin{equation}\label{45}
\Delta_{e}  (l_{X} s)  = D_{\bar{e}, X} s+  l_{[ \bar{e}, X]_{\nabla}}s. 
\end{equation}
On the other hand, from  the second relation in 
(\ref{e-unit-cond})
\begin{equation}
\Delta_{e}  (l_{X} s) = D_{\bar{e}, \bar{e}} (l_{X}s) = l_{X} (D_{\bar{e}, \bar{e}} s) + l_{\mathcal L_{\bar{e}} X}s=
D_{X, \bar{e}}s + l_{\mathcal L_{\bar{e}} X}s
\end{equation}
where in the second equality  we used  (\ref{integr-1})  and $\mathcal L_{X} (*)(\bar{e}, \bar{e}) = \mathcal L_{\bar{e}}X$,
while the third equaliy follows from Lemma \ref{lem-asoc}, 
by letting $Z:= X$ and $X = Y = \bar{e}$ in 
\begin{equation}
l_{Z} (D_{X, Y} s) + D_{X*Y, Z} s =l_{X} (D_{Z, Y} s) + D_{Z*Y, X} s.
\end{equation}
Using that 
$\mathcal L_{X} \bar{e} + [  \bar{e} , X]_{\nabla} $ equals  $2\nabla_{X} \bar{e}  +  T^{\nabla}_{\bar{e}, X}$
we obtain our claim.
\end{proof}

  \begin{lem}\label{duality-F} Assume that $\bar{\circ}$ is associative, with unit  field $e^{*}.$ Then $(E, \bar{\circ}, e^{*})$ is an $F$-manifold if and only if, for any
  $X, Y, Z, V\in {\mathfrak X}(M)$, 
  \begin{align}
 \nonumber&  [ \mathcal L_{Z} (X*Y) , V ]_{\nabla}  +  [\mathcal L_{V} (X*Y) , Z ]_{\nabla}  +\mathcal L_{ [X,Y]_{\nabla} } (*) (Z, V) \\
 \nonumber& - \mathcal L_{ [Z,V]_{\nabla} } (*) (X, Y) - [ X, \mathcal L_{Y} (*) (Z, V) ]_{\nabla}   - [ Y , \mathcal L_{X} (*) (Z, V)]_{\nabla} \\
  \nonumber& +\mathcal L_{Y}(*) (\mathcal L_{X} V, Z)  +\mathcal L_{Y}(*) (\mathcal L_{X} Z, V) 
+\mathcal L_{X}(*) (\mathcal L_{Y} V, Z)  +\mathcal L_{X}(*) (\mathcal L_{Y} Z, V)\\
\nonumber& + X* \left( [\mathcal L_{Y} V  , Z]_{\nabla}   + [ \mathcal L_{Y} Z  , V]_{\nabla}  \right)\\
\label{integr-duality} & + Y* \left( [ \mathcal L_{X} V  , Z]_{\nabla}   + [ \mathcal L_{X} Z  , V]_{\nabla} \right) 
  \end{align}
 belongs to the kernel of the map
  $$
  TM\ni X \rightarrow l_{X} \in \mathrm{End}\, (E).
  $$
  \end{lem}
  
  \begin{proof} We apply Lemma  \ref{prop-integr}.  Relation (\ref{integr-1}) for $(D^{*}, l^{*}, r^{*})$ follows from relation 
  (\ref{integr-1}) for $(D, l, r)$,  using the definition of $(D^{*}, l^{*}, r^{*})$  and $l_{X} ( l_{Y}s)= l_{X*Y} s$ for any $X$, $Y$. 
 We now consider  relation (\ref{F-man-lin}) for $(D^{*}, l^{*}, r^{*}).$  For any 
  $X, Y, Z, V\in {\mathfrak X}(M)$, 
  $s\in \Gamma (E^{*})$ and $\mu \in \Gamma (E^{*})$,  we define 
  \begin{align}
\nonumber T_{1}(X, Y, Z, V, \mu , s) & :=  (D^{*}_{\mathcal L_{X*Y} Z, V} \mu ) s  -   l^{*}_{X} ( D^{*}_{\mathcal L_{Y} V, Z} \mu )s -  l^{*}_{Y} (D^{*}_{\mathcal L_{X} V, Z} \mu )s,\\
\nonumber T_{2} (X, Y, Z, V,\mu , s)& :=  (D^{*}_{\mathcal L_{Y} (*) (Z, V) , X} \mu )s , \\
 \label{t} T_{3}(X, Y, Z, V, \mu , s) & :=  \left( D^{*}_{X, Y}  D^{*}_{Z, V}\mu   \right) s.
 \end{align}
 With this notation,  relation (\ref{F-man-lin}) for $(D^{*}, l^{*}, r^{*})$ becomes
 \begin{align}
\nonumber& T_{1}( X, Y, Z, V, \mu , s)  + T_{1} (X, Y, V, Z, \mu , s) + T_{2} (X, Y, Z, V, \mu , s) + T_{2} ( Y, X, Z, V, \mu , s) \\
\label{rel-t}& + T_{3} ( X, Y, Z, V, \mu , s) - T_{3} (Z, V, X, Y, \mu , s) =0. 
 \end{align} 
 From the definition of $D^{*}$ and $l^{*}$,  
 \begin{align}
 \nonumber& T_{1}(X, Y, Z, V, \mu , s) = ( \mathcal L_{X*Y} Z) ( \mu  l_{V}s)  - V (\mu l_{\mathcal L_{Z} (X* Y)} s ) 
 + (\mathcal L_{Z} (X*Y) * V ) (\mu s)\\
\nonumber&  - ( \mathcal L_{Y} V) ( \mu l_{X*Z} s) - Z(\mu l_{X *(\mathcal L_{Y} V)}s)  + ( (\mathcal L_{Y} V)* Z )  (  \mu l_{X}s) \\
\nonumber& - ( \mathcal L_{X} V)( \mu l_{Y*Z} s) - Z ( \mu l_{Y*( \mathcal L_{X} V)}s)  +  ((  \mathcal L_{X} V)* Z) (\mu l_{Y} s)\\
\nonumber& + \mu  l(s,  [\mathcal L_{Z} (X*Y), V]_{\nabla} + X* ( [ \mathcal L_{Y} V,  Z]_{\nabla} +Y* ( [ \mathcal L_{X} V,  Z]_{\nabla}) \\
 \nonumber&  -  \mu (  D_{\mathcal L_{X*Y} Z, V} s- D_{\mathcal L_{X}V, Z} (l_{Y}s) -  D_{\mathcal L_{Y} V, Z} (l_{X}s)),
  \end{align}
 where 
 $\mu s:= \mu (s)$ 
 and $l(s)(X) := l(s, X)$ for any  $s\in \Gamma (E)$,  $X\in {\mathfrak X}(M).$
By similar computations we obtain the other terms 
from  (\ref{rel-t}). Using  Lemmas \ref{lem-asoc} and \ref{prop-integr}  applied to $\circ$
we obtain, after cancellations,  that the left hand side of 
(\ref{rel-t}) equals $\mu (l_{\mathcal F} s)$, where $\mathcal F$ is given by (\ref{integr-duality}). 
\end{proof}

The next corollary concludes the proof of Theorem \ref{short}.

\begin{cor}\label{con-duality} Assume that $(M, *, \bar{e}, \nabla )$ is a flat $F$-manifold. Then 
$(E^{*}, \bar{\circ}, e^{*} )$ is a (FWL) $F$-manifold.
\end{cor}
\begin{proof}  Relations (\ref{asoc}) and (\ref{unit-duality}) follow from $T^{\nabla} =0$, $\nabla  \bar{e} =0$ and the symmetry of
$\nabla (*).$ We now prove that   expression (\ref{integr-duality}) vanishes, for any $X, Y, Z, V\in  {\mathfrak X}(M).$ 
Without loss of generality  we assume that $X, Y, Z, V$ are $\nabla$-parallel. 
As $T^{\nabla} =0$,  $X, Y, Z, V$ pairwise commute.
We need to show that 
$$
\nabla_{V} (\mathcal L_{Z} (X*Y)) +\nabla_{Z} (\mathcal L_{V} (X*Y))  = \nabla_{X} (\mathcal L_{Y} (Z*V)) + \nabla_{Y} (\mathcal L_{X} (Z*V)),
$$
or, using $T^{\nabla}=0$ again and $X, Y, Z, V\in {\mathfrak X}_{\nabla} (M)$, 
\begin{equation}\label{prime}
\nabla_{V}\nabla_{Z} (X*Y) +\nabla_{Z}\nabla_{V} (X*Y)  = \nabla_{X} \nabla_{Y} (Z*V) +\nabla_{Y} \nabla_{X} (Z*V).
\end{equation}
Since $R^{\nabla} =0$ and $\mathcal L_{V} Z =0$, 
the left hand side of (\ref{prime})  equals 
$$
2\nabla_{V} \nabla_{Z} (X*Y) = 2\nabla_{V} (\nabla_{Z}(*) (X, Y)) = 2\nabla_{V} (\nabla_{X} (*) (Z, Y) =2\nabla_{V}\nabla_{X} (Y*Z)
$$
where in the first equality we used $\nabla X  = \nabla Y =0$ and in the second the symmetry of $\nabla (*).$  Relation
(\ref{prime}) 
becomes $\nabla_{V}\nabla_{X} (Y*Z) = \nabla_{X}\nabla_{V} (Y*Z)$ and holds since 
$R^{\nabla} =0$ and $\mathcal L_{V} X =0.$ 
\end{proof}

\section{Enrichment of duality}\label{duality-enrich}
In this section we prove Theorem \ref{duality-euler-thm}.
Let   $\pi : E \rightarrow M$ be a vector bundle with base 
$(M, *, \bar{e},\nabla )$ a  flat $F$-manifold, 
$(E, \circ , e )$ a  FWL $F$-manifold over $(M, *,  \bar{e} )$   and $(E^{*}, \bar{\circ} , e^{*} ):= I^{\nabla}( E, \circ , e)$ 
(see  Theorem \ref{short}). Let $(D, l , *)$ be the components of $\circ .$

\subsection{Duality with  Euler fields}\label{duality-euler-sect}

The  claim on FWL $F$-manifolds with Euler fields in Theorem \ref{duality-euler-thm}  ii) follows from Theorem \ref{short} and  the next proposition. 

\begin{prop}\label{euler-duality}  The dual
$\mathcal E^{*}$ of a  FWL  Euler field  $\mathcal E$ on  $(E, \circ , e)$  is a (FWL)  Euler  field on $(E^{*}, \bar{\circ} , e^{*} )$ if and only if 
\begin{equation}\label{euler-duality-cond} 
l(s) ((\nabla^{2} \bar{\mathcal E} )_{X, Y} + (\nabla^{2} \bar{\mathcal E} )_{Y, X}) =0,
\end{equation}
for any  $s\in \Gamma (E)$ and $X, Y\in {\mathfrak X}(M)$, where $\bar{\mathcal E} := \pi_{*} \mathcal E .$
In particular, if $(M, *, \bar{e}, \nabla ,\bar{\mathcal E })$ is a flat $F$-manifold with Euler field, then $\mathcal E^{*}$ is 
an Euler field.
\end{prop}

\begin{proof}
We apply  Proposition \ref{euler-cor}.  As $\mathcal E^{*}$ projects to $\bar{\mathcal E}$ and the latter   is an Euler  field on $(M, *, \bar{e})$,
it remains to consider  the relations 
 \begin{align}
\nonumber& [ \Delta_{\mathcal E^{*}}, l^{*}_{X}]\mu  = l^{*}_{X}  \mu  + l^{*}_{\mathcal L_{\bar{\mathcal E}}X}\mu \\
\label{euler-cond-dual}& [ \Delta_{\mathcal E^{*}}, D^{*}_{X, Y} ] \mu   = D^{*}_{\mathcal L_{\bar{\mathcal E}}X, Y} \mu  + D^{*}_{\mathcal L_{\bar{\mathcal E}}Y, X} \mu 
+ D^{*}_{Y, Z} \mu ,
\end{align}
for any $X, Y\in { \mathfrak X}(M)$ and 
$\mu \in \Gamma (E^{*})$,  where  $(D^{*}, l^{*}, r^{*} )$ are the components of $\bar{\circ}$.
The first relation in (\ref{euler-cond-dual}) is a consequence of
the first relation in  (\ref{euler-cond}) (apply both sides   to $s\in \Gamma (E)$ and use the definition of $l^{*}$ and $\Delta_{\mathcal E^{*}}$
in terms of $l$ and $\Delta_{\mathcal E}$). 
We now consider the  second relation  in (\ref{euler-cond-dual}). Applying it  to $s\in \Gamma (E)$  and using
that $\bar{\mathcal E}$ is an Euler field on $(M, *, \bar{e})$, 
the first relation  in (\ref{euler-cond}) for $ [\Delta_{\mathcal E}, l_{X}]$, 
$ [\Delta_{\mathcal E}, l_{Y}]$ 
and $ [\Delta_{\mathcal E}, l_{[X, Y]_{\nabla} }]$,
and the second relation in (\ref{euler-cond})  for $[\Delta_{\mathcal E}, D_{X, Y} ]$,  we obtain that 
the  second relation in (\ref{euler-cond-dual}) is equivalent to
\begin{equation}\label{euler-duality-cond-1}
 l(s ) ((\mathcal L_{\bar{\mathcal E}} \nabla )_{X} Y + ( \mathcal L_{\bar{\mathcal E}}  \nabla )_{Y} X) =0,\ \forall  s\in \Gamma (E),\ X, Y\in {\mathfrak X}(M).
\end{equation}
Since $\nabla$ is flat and torsion-free,  $ (\mathcal L_{\bar{\mathcal E}} \nabla )_{X} Y = ( \nabla^{2} \bar{\mathcal E} )_{X, Y}$.
\end{proof}

\subsection{Duality with compatible connections}\label{duality-conn-sect}
Theorem \ref{duality-euler-thm}  i)
and the claim on FWL flat $F$-manifolds with Euler fields from Theorem \ref{duality-euler-thm} ii) 
follow from Corollaries \ref{flat-dual} and \ref{flat-euler-dual} below. 
Let  $\nabla^{fwl}$  be a FWL torsion-free connection on $E$, whose basic component $\nabla^{M}$ coincides with  
$\nabla .$ Let $\nabla^{E}$ and $\mathcal S$ be  the remaining components of $\nabla^{fwl} .$

\begin{defn}\label{dual-connection}  The {\cmssl dual  (FWL) connection} of  
$\nabla^{fwl}$ is the FWL torsion-free connection  $(\nabla^{fwl})^{*}$ 
on $E^{*}$ with components $(\nabla , \nabla^{E^{*}}, \mathcal S^{*})$, where 
$\nabla^{E^{*}}$ is the connection on the vector bundle $E^{*}$ induced by the connection $\nabla^{E}$ on the vector bundle $E$
and $\mathcal S^{*} \in \Gamma ( S^{2} T^{*}M\otimes \mathrm{End}\, E^{*})$ is defined by
\begin{equation}
(\mathcal S^{*}_{X, Y} \mu )s = - \mu  (\mathcal S_{X, Y} s),\ \forall \mu \in \Gamma (E^{*}),\ s\in \Gamma (E). 
\end{equation}
\end{defn}

\begin{prop}\label{compat-dual} Assume that  $\nabla^{fwl}$ is compatible on $(E, \circ , e)$. Then  $(\nabla^{fwl})^{*}$ is compatible on
$(E^{*}, \bar{\circ}, e^{*})$ if and only if, for any $X, Y, Z\in {\mathfrak X}(M)$,
\begin{equation}\label{comutatori} 
[ l_{X}, \mathcal S_{Y, Z} ] - [ l_{Y}, \mathcal S_{X, Z} ] = \frac{1}{2} \left( [ l_{X}, R^{\nabla^{E}}_{Y, Z} ] -  [l_{Y}, R^{\nabla^{E}}_{X, Z} ] \right)
+ [ l_{Z}, R^{\nabla^{E}}_{X, Y} ].
\end{equation}
\end{prop}

\begin{proof}   We apply Subsection \ref{duality-conn-sect}. 
Relation (\ref{con-2}) applied to $(\nabla^{fwl} )^{*}$ 
is 
\begin{equation}
\nabla^{E^{*}}_{Y*Z} \mu + D^{*}_{Y, Z} \mu - l_{Y}^{*} (\nabla^{E^{*}}_{Z}\mu  )-\nabla^{E^{*}}_{Y} (l_{Z}\mu ) + l_{\nabla_{Y}Z}\mu =0,
\end{equation}
for any $\mu \in \Gamma (E^{*})$ and $Y, Z\in {\mathfrak X}(M)$.
Applied to $s\in \Gamma (E)$, it becomes
\begin{equation}\label{e-d-0}
\nabla^{E}_{Y*Z} s+ D_{Y, Z} s+ l_{\nabla_{Z}Y} s  -\nabla^{E}_{Z} ( l_{Y} s) - l_{Z} (\nabla^{E}_{Y} s) =0.
 \end{equation}
From relation (\ref{con-2}) applied to $\nabla^{fwl}$,  the  sum of the first two terms in  the left hand side of (\ref{e-d-0}) equals
$l_{Y} \nabla^{E}_{Z} s + \nabla^{E}_{Y} (l_{Z} s) - l_{\nabla_{Y}Z}s$. We deduce that (\ref{e-d-0})  coincides with
(\ref{con-3}) applied  to  $\nabla^{fwl}$ and is therefore satisfied.  Similarly,  (\ref{con-3})  and  the first relation in (\ref{nabla-fwl-e}) 
applied  to  $(\nabla^{fwl})^{*}$ follow from the same relations  applied to $\nabla^{fwl}$. The same is true for the second
relation in  (\ref{nabla-fwl-e})  applied to $\nabla^{fwl}$ and  $(\nabla^{fwl})^{*}$, by using 
\begin{equation}\label{curvature}
(R^{\nabla^{E^{*}}}_{X, Y} \mu )s = - \mu ( R^{\nabla^{E}}_{X, Y} s).
\end{equation}
 Using (\ref{curvature}) again,   we obtain that 
$(\nabla^{fwl})^{*}$ satisfies (\ref{con-1-alt})
if and only if
\begin{align}
\nonumber&\mathcal S_{X, Y*Z} s- \mathcal S_{Y, X*Z} s-  \mathcal S_{X, Z}(l_{Y} s)+ 
\mathcal S_{Y, Z}(l_{X}s)\\
\nonumber& = \frac{1}{2} \left( R^{\nabla^{E}}_{X, Y*Z}s - R^{\nabla^{E}}_{Y, X*Z}s -
R^{\nabla^{E}}_{X, Z}(l_{Y} s) + R^{\nabla^{E}}_{Y, Z}(l_{X}s)\right) - l_{Z} (R^{\nabla^{E}}_{X, Y} s).
\end{align}
Combining the above relation with  (\ref{con-1-alt}) applied to $\nabla^{fwl}$ we obtain (\ref{comutatori}).
\end{proof}

\begin{cor}\label{flat-dual}  If $ (E, \circ , e, \nabla^{fwl})$ is a flat $F$-manifold  then $ (E^{*}, \bar{\circ} , e^{*},  (\nabla^{fwl})^{*})$ is a flat $F$-manifold. 
\end{cor}

\begin{proof} Since $\nabla^{fwl}$ is flat, $\nabla^{E}$ is flat and $\mathcal S =0$ (from  Corollary \ref{flat-fwl}). 
Then  $(\nabla^{fwl})^{*}$ is also flat. 
From  (\ref{comutatori}),  $(\nabla^{fwl})^{*}$ is compatible on $ (E^{*}, \bar{\circ} , e^{*}).$ 
\end{proof}

\begin{lem}\label{twice} Let $\mathcal E \in {\mathfrak X}_{\mathrm{fwl}}(E)$ (not necessarily an Euler field) such that 
$(\nabla^{fwl})^{2} \mathcal E =0.$  Then $((\nabla^{fwl})^{*})^{2} \mathcal E^{*} =0$ if and only if, for any $X, Y
\in {\mathfrak X}(M)$,
\begin{equation}\label{E-der-req}
\{ \Delta_{\mathcal E} -\nabla_{\bar{\mathcal E}}, \mathcal S_{X, Y} -\frac{1}{2} R^{\nabla^{E}}_{X, Y} \} =0,\ \forall X, Y\in {\mathfrak X}(M),
\end{equation} 
 where $\{ A, B\} := AB + BA$ denotes the anti-commutator of $A, B\in \Gamma (\mathrm{End}\, E)$.
\end{lem}

\begin{proof} Apply Corollary \ref{cor-E}  and notice that (\ref{cor-E-rel}) on the dual side gives (\ref{E-der-req}).\end{proof}

\begin{cor} \label{flat-euler-dual} Assume that $(E, \circ , e,  \nabla^{fwl}, \mathcal E )$ is a flat $F$-manifold with Euler field. Then 
$(E^{*}, \bar{\circ} , e^{*},  (\nabla^{fwl} )^{*}, \mathcal E^{*})$ is a flat $F$-manifold with Euler field.
\end{cor}

\begin{proof}  From Corollary \ref{flat-dual},   
$(E, \bar{\circ} , e^{*},  (\nabla^{fwl})^{*})$ is a flat $F$-manifold. From Proposition \ref{euler-duality},
$\mathcal E^{*}$ is an Euler field. 
From  Lemma \ref{twice}, $((\nabla^{fwl})^{*})^{2} \mathcal E^{*} =0$ as 
$R^{\nabla^{E}}=0$ and $\mathcal S =0$.
\end{proof}

\section{Examples  in small dimensions}\label{examples-sect}

 \begin{exa}\label{ex-1}{\rm  This example is inspired from  Theorem 4.7  (b) 
of \cite{h-book}.  Let $\pi : E \rightarrow \mathbb{R}$ be a  trivial   line bundle  and 
$(x, \xi )$ coordinates on $E$ 
determined by a trivialization 
$s\in \Gamma (E)$,  with $x$ the base  and $\xi$ the fiber coordinate.
Define $(D, l, * )$ by 
$D_{\partial x, \partial x} (f s) = f^{\prime}(x) s$
for any $f\in C^{\infty}(\mathbb{R})$,  $l_{\partial x} s = s$
and $\partial x * \partial x = \partial x$. Then $(D, l, *)$ are components of a FWL  multiplication  $\circ$ on $E$, with  FWL unit field $e= \partial x$ and 
$\partial \xi \circ \partial \xi =0$. We obtain a FWL  $F$-manifold  $(E, \circ , e)$. Any FWL   Euler field  is of the form
 $\mathcal E  = ( x + c) \partial x   +  a\xi  \partial \xi $
 where $a, c\in \mathbb{R}$.  The  connection $\nabla$    defined by 
  $\nabla \partial x =0$  makes $(\mathbb{R}, *, \partial x)$ a flat $F$-manifold.
 The dual   multiplication $\bar{\circ}$  and $\mathcal E^{*}$ have the same form  as $\circ$ 
 and $\mathcal E$ (in local coordinates  of $E^{*}$ determined by  $\mu \in \Gamma (E^{*})$, $\mu (s) =1$) 
 with $a$ replaced by $-a$ and $c$ unchanged. In particular,   $\mathcal E^{*}$ is Euler  on $(E^{*}, \bar{\circ}, e^{*})$. }
 \end{exa}

\begin{exa}{\rm  This example is inspired  from Theorem 5.4 of \cite{h-book}.
Let $\pi : E \rightarrow \mathbb{R}^{2}$ be a trivial line  bundle,
$s\in \Gamma (E)$ a trivialization  and $(x_{1}, x_{2}, \xi )$ the corresponding coordinates on $E$, with $(x_{1}, x_{2})$  the base and $\xi$ the fiber coordinate.  Let $h\in C^{\infty}(  \mathbb{R})$  
and define $(D, l , *)$ 
by $l_{\partial x_{1}} s =s$, $l_{\partial x_{2}} s =0$, 
$\partial x_{1} * \partial x_{i} = \partial x_{i}$ ($i=1,2$), $\partial x_{2} * \partial x_{2} =0$,
$D s = h(x_{2}) dx_{2} \otimes dx_{2} \otimes s$
(extended to $\Gamma (E)$ using
(\ref{comp-cond})).  
Then  $(D, l, * )$ are components of  a FWL multiplication $\circ$ on $E$, with FWL unit field
$ e = \partial x_{1}$ and 
\begin{equation}\label{m-d}
\partial x_{2} \circ  \partial x_{2} = \xi h(x_{2}) \partial \xi ,\ \partial x_{2} \circ \partial \xi = \partial \xi \circ \partial \xi =0.
\end{equation}
We obtain a FWL $F$-manifold $(M, \circ , e)$ over the $F$-manifold from Example \ref{ex-1}.
Any vector field of the form
\begin{equation}\label{e-d}
\mathcal E = ( x_{1}  + c) \partial x_{1}  + \epsilon  (x_{2} )\partial x_{2}   +  \xi   \eta ( x_{2})  \partial  \xi
\end{equation}
where $c\in \mathbb{R}$,  $\epsilon$, $\eta$ depend only on $x_{2}$ and $\epsilon$ satisfies
\begin{equation}\label{E1}
\epsilon (x_{2}) h^{\prime} (x_{2}) + 2\epsilon^{\prime} (x_{2}) h(x_{2}) = h(x_{2})
\end{equation}
is a FWL Euler field. Any connection $\nabla$ with $\nabla \partial x_{1} =0$,  $\nabla_{\partial x_{1}} \partial x_{2} = 0$ and
$\nabla_{\partial x_{2}} \partial x_{2} =\lambda \partial x_{1} +\beta \partial x_{2}$ where $\lambda$, $\beta$ depend only on $x_{2}$, makes 
$(\mathbb{R}^{2}, *, \partial x_{1}, \nabla )$ a flat $F$-manifold. Remark that  $\nabla^{2} \bar{\mathcal E} =0$ if and only if
\begin{equation}\label{E2}
(\epsilon^{\prime} +\epsilon \beta )^{\prime} =0,\ 2\epsilon^{\prime} \lambda + \epsilon \lambda^{\prime} =\lambda .
\end{equation}
The dual  multiplication $\bar{\circ}$ 
and $\mathcal E^{*}$ have  the same form as (\ref{m-d})
and (\ref{e-d})  (in local coordinates of $E^{*}$ determined by $\mu\in \Gamma (E^{*})$, $\mu (s) =1$), with $h$ replaced by $ - (h+2\lambda ) $,
 $\eta$ by $-\eta$ and $c$, $\epsilon$ unchanged. 
When (\ref{E2}) are satisfied,  $\mathcal E^{*}$  is Euler on $(E^{*}, \bar{\circ}, e)$.  This can also be checked directly: 
the second relation in (\ref{E2}) implies that 
if $(\epsilon , h)$ satisfy (\ref{E1}), then also $(\epsilon, - (h+2\lambda ))$ satisfy (\ref{E1}).}
\end{exa}

\section{The tangent prolongation}\label{tangent-sect}

\subsection{Proof of Proposition \ref{integr-comp}}
The identity 
$$
\mathcal L_{fZ} ( * ) (X, Y) = f \mathcal L_{Z} ( * ) (X, Y) + X(f) Z *  Y + Y(f) Z * X - (X * Y)(f) Z,
$$
for any $X, Y, Z\in {\mathfrak X}(M)$ and $f\in C^{\infty}(M)$  implies that 
$(D, l, * )$, defined by (\ref{prolong}), are the components of a FWL commutative multiplication $*^{T}$ on $TM$
(see Proposition \ref{comp-prop}).
Using Lemmas  \ref{lem-asoc} and \ref{lem-unit},  one can verify  that $*^{T}$ is associative with unit  field ${e}^{T}.$ 
E.g., to check that ${e}^{T}$ is a unit  field 
for $*^{T}$ we need  to verify that 
$\pi_{*} ({e}^{T})$ is a unit  field for $*$,  
$l_{\pi_{*} ({e}^{T}) } (X) = X$ (which hold  as $\pi_{*} ({e}^{T}) =e$)
and  $D_{X, e} Y = l_{X} \Delta_{{e}^{T}} Y$, 
for any $X, Y\in {\mathfrak X}(M).$ 
The latter relation is equivalent to
$\mathcal L_{Y} (* ) (X,  e ) = (\mathcal L_{e} Y)*   X$ and  is satisfied.
It remains to check that $(D, l, *)$ satisfy 
(\ref{integr-1}) and (\ref{F-man-lin}). Relation (\ref{integr-1}) reduces to  (\ref{integr-F-man}) 
and   (\ref{F-man-lin}) to (\ref{E}) below.  The next lemma concludes the proof of 
Proposition \ref{integr-comp}.

\begin{lem}\label{prolong-aux}  On any $F$-manifold $(M, *, e)$, 
\begin{align}
\nonumber& 2 \mathrm{sym}_{(Z, V)} \mathcal L_{W} (*) (\mathcal  L_{X*Y} Z, V) +2 \mathrm{sym}_{(X, Y)} \mathcal L_{W} 
(*) (\mathcal L_{Y} (*) (Z, V), X)\\
\nonumber&  - 4\mathrm{sym}_{(X, Y)}\mathrm{sym}_{(Z, V)} X* \mathcal L_{W} (*) (\mathcal L_{Y} V, Z) \\
\label{E}& + \mathcal L_{\mathcal L_{W} (*) (Z, V)} (*)(X, Y)  -  \mathcal L_{\mathcal L_{W} (*) (X, Y)} (*)(Z, V)
=0,
\end{align}
for any  $X, Y, Z, V, W\in {\mathfrak X}(M)$,
where $\mathrm{sym}_{(X, Y)} F(X, Y, Z, V, W)$   denotes the symmetric  
part in   $X$,  $Y$ of an expression $F(X, Y, Z, V, W)$  depending on $X, Y, Z, V,W$  (similarly for $\mathrm{sym}_{(Z, V)} F(X, Y, Z, V,W)$).
\end{lem}
\begin{proof}  Let $\mathcal T$ be  the first line in (\ref{E}). 
Expanding 
$$
\mathcal L_{W} (*) (\mathcal L_{X*Y}Z, V) = \mathcal L_{W} ((\mathcal L_{X*Y}Z)* V) - (\mathcal L_{W} \mathcal L_{X*Y} Z)* V
-(\mathcal L_{X*Y} Z)* \mathcal L_{W} V
$$ 
and similarly 
\begin{align}
\nonumber& \mathcal L_{W} (*) (\mathcal L_{Y} (*) (Z, V), X)= \mathcal L_{W} ( (\mathcal L_{Y} (*) (Z, V) * X)
- \mathcal L_{W} (\mathcal L_{Y} (*) (Z, V)) * X\\
\nonumber&  - \mathcal L_{Y}(*) (Z, V) * \mathcal L_{W} X
\end{align}
and using
$$
(\mathcal L_{X*Y} Z) *V + (\mathcal L_{X*Y} V)*Z + \mathcal L_{Y}(*) (Z, V) * X +\mathcal L_{X}(*)(Z, V) * Y = {\mathcal L}_{X*Y} (Z*V)
$$
for the terms inside $\mathcal L_{W}$, we obtain 
\begin{align*}
\nonumber&\mathcal T = 
\mathcal L_{W} (\mathcal L_{X*Y} (Z*V) ) -\mathcal L_{W} (\mathcal L_{X*Y}Z) * V -( \mathcal L_{X*Y}Z)*\mathcal L_{W}V\\
\nonumber&- \mathcal L_{W} (\mathcal L_{X*Y} V) * Z -(\mathcal L_{X*Y} V)* \mathcal L_{W} Z -\mathcal L_{W} (\mathcal L_{Y}(*) (Z, V)) * X\\
\nonumber& -\mathcal L_{Y} (*) (Z, V) * \mathcal L_{W}X -\mathcal L_{W} (\mathcal L_{X}(*) (Z, V)) * Y -\mathcal L_{X}(*) (Z, V) *
\mathcal L_{W}Y.
\end{align*}
Using the Jacobi identity followed by  (\ref{integr-F-man}) we rewrite the first term in $\mathcal T$: 
\begin{align}
\nonumber& \mathcal L_{W} \mathcal L_{X*Y}  (Z*V) = \mathcal L_{ \mathcal L_{W} (X*Y)} (Z*V) +\mathcal L_{\mathcal L_{Z*V}W} (X*Y)\\
\nonumber& = \mathcal L_{\mathcal L_{W}(*)(X, Y)} (*)(Z, V) +(\mathcal L_{W}X)* \mathcal L_{Y} (*) (Z, V) + Y *\mathcal L_{\mathcal L_{W} X}(*)
(Z, V)\\
\nonumber& + X* \mathcal L_{\mathcal L_{W}Y} (*) (Z, V) +(\mathcal L_{W} Y)*\mathcal L_{X}(*) (Z, V)
+(\mathcal L_{\mathcal L_{W}(X*Y)}Z)* V\\
\nonumber&  +(\mathcal L_{\mathcal L_{W}(X*Y)}V)* Z -\mathcal L_{\mathcal L_{W} (*) (Z, V)}(*) (X, Y) - (\mathcal L_{W} Z) * \mathcal L_{V}(*) (X, Y)\\
\nonumber& 
- V* \mathcal L_{\mathcal L_{W}Z}(*) (X, Y) - Z* \mathcal L_{\mathcal L_{W}V}(*) (X, Y)  - (\mathcal L_{W} V)* \mathcal L_{Z}(*) (X, Y) \\
\nonumber& +(\mathcal L_{\mathcal L_{Z*V}W} X)*Y + (\mathcal L_{\mathcal L_{Z*V}W} Y)*X.
\end{align}
Replacing the above expression in $\mathcal T$ we obtain
\begin{align}
\nonumber& \mathcal T = \mathcal L_{\mathcal L_{W} (*) (X, Y)} (*) (Z, V) -\mathcal L_{\mathcal L_{W}(*) (Z, V)} (*) (X, Y) \\
\nonumber& 
- (\mathcal L_{W} V)* \left( \mathcal L_{Z}(*) (X, Y) +\mathcal L_{X*Y} Z\right) - (\mathcal L_{W} Z)* \left( \mathcal L_{V} (*) (X, Y) +\mathcal L_{X*Y} V\right)\\
\nonumber& +Y *\left( \mathcal L_{\mathcal L_{W} X}(*) (Z, V) +\mathcal L_{\mathcal L_{Z*V}W} X - \mathcal L_{W} ( \mathcal L_{X} (*) (Z, V))\right)\\
\nonumber& + X* \left(  \mathcal L_{\mathcal L_{W} Y}(*) (Z, V) +\mathcal L_{\mathcal L_{Z*V}W} Y - \mathcal L_{W} ( \mathcal L_{Y} (*) (Z, V))\right)\\
\nonumber& + V * \left( \mathcal L_{\mathcal L_{W} (X*Y)} Z -\mathcal L_{\mathcal L_{W} Z} (*) (X, Y) -\mathcal L_{W} ( \mathcal L_{X*Y} Z)\right)\\
\label{E-rel}& + Z * \left( \mathcal L_{\mathcal L_{W} (X*Y)} V -\mathcal L_{\mathcal L_{W} V} (*) (X, Y) -\mathcal L_{W} ( \mathcal L_{X*Y} V)\right) .
\end{align}
In the last two lines 
we  replace  $\mathcal L_{\mathcal L_{W} (X*Y)} Z -\mathcal L_{W} ( \mathcal L_{X*Y} Z)$ by
$-\mathcal L_{\mathcal L_{Z} W} (X*Y)$ and the same with $V$ and $Z$ interchanged and we obtain a new expression of $\mathcal T$.
Replacing it in  (\ref{E}), the latter  becomes 
\begin{align}
\nonumber& \mathrm{sym}_{(X, Y)} X* \left( 
(\mathcal L_{W} V)*(\mathcal L_{Z} Y) +  (\mathcal L_{W} Z)*(\mathcal L_{V} Y)
+\mathcal L_{\mathcal L_{W} Y} (*) (Z, V)+\mathcal L_{\mathcal L_{Z*V}W}Y\right)\\
\nonumber&- \mathrm{sym}_{(X, Y)} X* \left( \mathcal L_{W} (\mathcal L_{Y}(*) (Z, V)) 
+\mathcal L_{W}(*) (\mathcal L_{Y}V, Z) +\mathcal L_{W}(*) (\mathcal L_{Y}Z, V) \right) \\
\nonumber& - \mathrm{sym}_{(Z, V)} Z* \left( \mathcal L_{\mathcal L_{V}W} (X*Y) +\mathcal L_{\mathcal L_{W} V}(*) (X, Y)\right)=0.
\end{align}
We now simplify the first two lines in the above relation, using 
\begin{align}
\nonumber&  (\mathcal L_{W} V)*(\mathcal L_{Z} Y) +  (\mathcal L_{W} Z)*(\mathcal L_{V} Y)
+\mathcal L_{\mathcal L_{W} Y} (*) (Z, V)+\mathcal L_{\mathcal L_{Z*V}W}Y\\
\nonumber&  - \mathcal L_{W} (\mathcal L_{Y}(*) (Z, V)) 
- \mathcal L_{W}(*) (\mathcal L_{Y}V, Z) - \mathcal L_{W}(*) (\mathcal L_{Y}Z, V)\\
\nonumber& = (\mathcal L_{Y} \mathcal L_{W} Z)*V +(\mathcal L_{Y} \mathcal L_{W} V)*Z,
\end{align}
which is obtained by expanding  $\mathcal L_{Y}(*)(Z, V) $ inside $\mathcal L_{W} (\mathcal L_{Y}(*) (Z, V)) $, 
using  the Jacobi identity for the term $\mathcal L_{W} (\mathcal L_{Y}(Z*V)) $ and  the identity 
$$
\mathcal L_{Y} (Z*V)= \mathcal L_{Y}(*)(Z, V)  + (\mathcal L_{Y} Z)*V + Z*(\mathcal L_{Y} V)
$$ 
for various arguments $Y, Z, V.$  Relation (\ref{E}) becomes 
\begin{align}
\nonumber& Z * \left( Y* \mathcal L_{X} \mathcal L_{W} V + X* \mathcal L_{Y} \mathcal L_{W} V -\mathcal L_{\mathcal L_{V} W} (X*Y) -
\mathcal L_{\mathcal L_{W} V}(*) (X, Y)\right)\\
\nonumber&  + V * \left( Y* \mathcal L_{X} \mathcal L_{W} Z + X* \mathcal L_{Y} \mathcal L_{W} Z -\mathcal L_{\mathcal L_{Z} W} (X*Y) -\mathcal L_{\mathcal L_{W} Z}(*) (X, Y)\right) =0
\end{align}
and  holds as both lines from  the left hand side  vanish. 
\end{proof}

\subsection{Proof of Proposition \ref{e-con} i)}

We  apply Proposition \ref{euler-cor}. The first relation in 
(\ref{euler-cond})  follows directly  from $\mathcal L_{\mathcal E} (*) = *$. The second relation in (\ref{euler-cond}) becomes
\begin{align}
\nonumber& \mathcal L_{\mathcal E} ( \mathcal L_{V}(* ) (Y, Z)) -\mathcal L_{\mathcal L_{\mathcal E}V} (*) (Y, Z) - \mathcal L_{V} (*)
 (\mathcal L_{\mathcal E} Y, Z) -\mathcal L_{V} (*)(Y, \mathcal L_{\mathcal E} Z)\\
 \label{euler-tg}& = \mathcal L_{V} (*) (Y, Z),
\end{align}
for any $V, Y, Z\in {\mathfrak X}(M).$ 
Proposition \ref{e-con} i) follows from the next lemma.

\begin{lem} Any Euler field $\mathcal E$ on an $F$-manifold  $(M, *, e)$ satisfies  (\ref{euler-tg}).
\end{lem}

\begin{proof} Expanding   all terms  of the form 
$\mathcal L_{X}(*) (Y, Z)$  from the left hand side of 
(\ref{euler-tg})  and noticing that 
\begin{align}
\nonumber& \mathcal L_{\mathcal E} ( \mathcal L_{V}(* ) (Y, Z))= \mathcal L_{\mathcal E} \mathcal L_{V} (Y*Z) - (\mathcal L_{V}Y)*Z
- Y* (\mathcal L_{V} Z)\\
\nonumber& -  (\mathcal L_{\mathcal E} \mathcal L_{V} Y)* Z - Y* \mathcal L_{\mathcal E} \mathcal L_{V} Z -(\mathcal L_{V} Y)
*( \mathcal L_{\mathcal E} Z) - (\mathcal L_{\mathcal E} Y ) *( \mathcal L_{V} Z)
\end{align}
we obtain that the left-hand side of (\ref{euler-tg}) equals 
\begin{align}
\nonumber L.H.S&= \mathcal L_{\mathcal E}\mathcal L_{V} (Y*Z) + Z* ( - \mathcal L_{V} Y -\mathcal L_{\mathcal E} \mathcal L_{V} Y +\mathcal L_{\mathcal L_{\mathcal E} V} Y +\mathcal L_{V}\mathcal L_{\mathcal E} Y)\\
\nonumber& + Y* ( - \mathcal L_{V} Z -\mathcal L_{\mathcal E} \mathcal L_{V} Z +\mathcal L_{\mathcal L_{\mathcal E}V} Z +\mathcal L_{V}\mathcal L_{\mathcal E} Z)\\
\label{lhs}& - \mathcal L_{\mathcal L_{\mathcal E}V} (Y*Z) - \mathcal L_{V} ((\mathcal L_{\mathcal E} Y) *Z) -\mathcal L_{V} 
(Y* \mathcal L_{\mathcal E}Z).
\end{align}
From the Jacobi identity, the  sum of the last three terms in  the first and second lines in the right hand side of (\ref{lhs})  vanish and we obtain that 
\begin{align}
\nonumber L.H.S& = \mathcal L_{\mathcal E} \mathcal L_{V} (Y*Z) - Z* \mathcal L_{V} Y  -  Y*  \mathcal L_{V} Z \\
\nonumber& - \mathcal L_{\mathcal L_{\mathcal E}V} (Y*Z) - \mathcal L_{V} ((\mathcal L_{\mathcal E} Y) *Z) -\mathcal L_{V} 
(Y* \mathcal L_{\mathcal E}Z).
\end{align}
 From the Jacobi identity again, 
$\mathcal L_{\mathcal E} \mathcal L_{V} (Y*Z) - \mathcal L_{\mathcal L_{\mathcal E}V} (Y*Z)= 
\mathcal L_{V} \mathcal L_{\mathcal E} (Y*Z)$
which implies that
\begin{align}
\nonumber L.H.S& = \mathcal L_{V} \mathcal L_{\mathcal E} (Y*Z) - Z* \mathcal L_{V} Y  -  Y*  \mathcal L_{V} Z \\
\nonumber& - \mathcal L_{V} ((\mathcal L_{\mathcal E} Y) *Z) -\mathcal L_{V} 
(Y* \mathcal L_{\mathcal E}Z)\\
\nonumber&  = \mathcal L_{V} (\mathcal L_{\mathcal E} (*) (Y, Z))  - Z* \mathcal L_{V} Y  -  Y*  \mathcal L_{V} Z \\
\nonumber& = \mathcal L_{V} (Y * Z)  - Z* \mathcal L_{V} Y  -  Y*  \mathcal L_{V} Z = \mathcal L_{V}(*) (Y, Z)
\end{align}
where in the last equality  we used that $\mathcal E$ is an Euler field. 
\end{proof}

\subsection{Proof of Proposition \ref{e-con} ii), iii)} 

We start with a general lemma. 

\begin{lem}\label{complete} Let  $\nabla$ be a  torsion-free connection on a manifold $M$. There is a  unique FWL torsion-free  connection  
$\nabla^{T}$ on $TM$, with
components $(\nabla , \nabla , \mathcal S )$, which satisfies 
\begin{equation}\label{T}
\nabla^{T} _{X^{T}} Y^{T} =( \nabla_{X}Y)^{T},\ \forall X, Y\in {\mathfrak X}(M).
\end{equation}
Its $\mathcal S$-component is given by 
\begin{equation}\label{mathcalG} 
\mathcal S_{X, Y} Z = R^{\nabla}_{Z, X} Y +\frac{1}{2} R^{\nabla}_{X, Y} Z.
\end{equation}
\end{lem}

\begin{proof}  Let $(x^{i})$ be local coordinates on $M$ and $(x^{i}, u^{j})$ the induced coordinates on $TM.$ 
Write
\begin{align}
 \nonumber& \nabla^{T}_{\partial x^{i} }( \partial x^{j} ) = \Gamma_{ij}^{k} (x) \partial x^{k}  + 
 \Gamma_{ij\beta}^{\alpha} u^{\beta} \partial u^{\alpha}\\
  \label{T-1}& \nabla^{T}_{\partial x^{i} } (\partial u^{j} ) = \nabla^{T}_{\partial u^{j}} 
 (\partial x^{i} )
  = \bar{\Gamma}_{ij}^{\gamma}(x) \partial u^{\gamma},
  \end{align} 
 where $\Gamma_{ij}^{k}$ (the Christoffel symbolds of $\nabla$ in coordinates $(x^{i})$),  $\Gamma_{ij\beta}^{\alpha} $ and  $\bar{\Gamma}_{i\beta}^{\gamma}$ are local functions on $M$.
 Recall that all other covariant derivatives of coordinate vector fields in  direction of coordinate vector fields vanish.
Let  $X = X^{i}\partial x^{i}$ and 
$Y = Y^{i} \partial x^{i}$ be vector fields on $M$.
Then  $X^{T} = X^{i} \partial x^{i} + u^{j} \frac{\partial X^{i}}{\partial x^{j}} 
\partial u^{i}$ and similarly for $Y^{T}$. 
Replacing $X^{T}$ and $Y^{T}$ in 
(\ref{T}) we obtain, using   (\ref{T-1}),  that $\Gamma_{ip\beta}^{\alpha} =\frac{\partial \Gamma^{\alpha}_{ip}}{\partial x^{\beta}}$. So
the  $\mathcal G$-component of $\nabla^{T}$ is given by 
$\mathcal G (\partial x^{\beta})
(\partial x^{i}, \partial x^{j}) =\frac{\partial \Gamma^{\alpha}_{ij}}{\partial x^{\beta}} \partial x^{\alpha}$
(as
$\mathcal L_{\partial u^{\beta}} (\nabla^{T}) (\partial x^{i}, \partial x^{j}) =\Gamma^{\alpha}_{ij\beta} \partial u^{\alpha}$). 
Using that $\nabla$ is torsion-free, 
 \begin{align}
\nonumber&  \mathcal G (\partial x^{\beta})
(\partial x^{i}, \partial x^{j}) =
\nabla_{\partial x^{\beta}}  (\Gamma^{\alpha}_{ij} \partial x^{\alpha}) - \Gamma^{\alpha}_{ij}
\nabla_{\partial x^{\beta} }  (\partial x^{\alpha})\\
\nonumber&=  \nabla_{\partial x^{\beta} } \nabla_{\partial x^{i} } (\partial x^{j})
- dx^{\alpha} ( \nabla_{\partial x^{i}} (\partial x^{j})) \nabla_{\partial x^{\alpha}} (\partial x^{\beta}) \\
\nonumber& =  R^{\nabla}_{\partial x^{\beta}, \partial x^{i}} (\partial x^{j}) 
+\nabla_{\partial x^{i}} \nabla_{\partial x^{\beta }} (\partial x^{j}) -
\nabla_{\nabla_{\partial x^{i}} ( \partial x^{j})} (\partial x^{\beta} )\\
\nonumber& = R^{\nabla}_{\partial x^{\beta}, \partial x^{i}} (\partial x^{j}) 
+\nabla^{2} _{\partial x^{i}, \partial x^{j}}
(\partial x^{\beta}).
\end{align}
From  (\ref{g-s-1}), we obtain  (\ref{mathcalG}) with $X, Y, Z$ coordinate vector fields.
\end{proof}

The next proposition implies  Proposition \ref{e-con} ii), iii).

\begin{prop}\label{compat-tg}
Let $(M, *, e)$ be an $F$-manifold and $\nabla$ a compatible connection.
The connection $\nabla^{T}$  is compatible  on $(TM, *^{T}, e^{T})$ if and only if 
\begin{align}
\nonumber& R^{\nabla}_{V, X} (Y*Z) - R^{\nabla}_{V, Y} (X*Z) + R^{\nabla}_{X, Y} (Z*V)\\
\label{con-n-comp}& + (R^{\nabla}_{V, Y} Z)* X - (R^{\nabla}_{V , X} Z)* Y - (R^{\nabla}_{X, Y} Z) * V =0,
\end{align}
for any $X, Y, Z, V\in {\mathfrak X}(M).$ If $\nabla$ is flat, then  so is $\nabla^{T}$.
If $(M, *, e , \nabla , \mathcal E)$ is a flat $F$-manifold with Euler field, then so is 
$(TM, * , e^{T} , \nabla^{T} , \mathcal E^{T})$.
\end{prop}

\begin{proof} 
The flatness of $\nabla^{T}$ (when $\nabla$ is flat) follows  from  Corollary \ref{flat-fwl} and  (\ref{mathcalG}).
From Proposition \ref{e-con} i), $\mathcal E^{T}$ is an Euler field. 
The basic component of $(\nabla^{T} )^{2} (\mathcal E^{T})\in
\mathcal T^{2,1}_{\mathrm{fwl}}(TM)$ is $\nabla^{2}\mathcal E $ and the remaining components  
are given by 
\begin{equation}
l^{(1)}(X, Y) = l^{(2)}(X, Y) = (\nabla^{2} \mathcal E )_{X,Y},\ D_{X, Y} Z = \mathcal L_{Z} (\nabla^{2}\mathcal E )_{X, Y}.
\end{equation}
We obtain that if
$\nabla^{2} \mathcal E =0$ then also $(\nabla^{T} )^{2} (\mathcal E^{T})=0.$ 
For  the first claim, we use 
Propositions  \ref{nabla-circ-sym},  \ref{mathcal-e} 
and Lemma \ref{mathcal-s}.  All relations from these statements are satisfied by
$(*^{T}, e^{T}, \nabla^{T})$ except relation (\ref{con-1-alt}) (or (\ref{con-1})). Using
(\ref{mathcalG}), we obtain that (\ref{con-1-alt}) 
reduces to  (\ref{con-n-comp}).
\end{proof}

\section{Further  prolongations of $F$-manifolds}\label{ctg-def} 

{\bf The cotangent and generalized prolongations.}  We define the {\cmssl cotangent prolongation}  
$(T^{*}M,  *^{T^{*}}, e^{T^{*} })$ 
of a flat $F$-manifold $(M, *, e, \nabla )$ as the image of the tangent prolongation 
$(TM,  *^{T} , e^{T})$   in the duality  $I^{\nabla}$ of Theorem \ref{short}.  
For the generalized   prolongation we need  the next lemma, which is a consequence of the previous sections (see also Remark \ref{iso-sum}): 

\begin{lem}\label{iso-sum-cor}   Let $\pi_{i} : E_{i} \rightarrow M$  be vector bundles ($i=1,2$). 
If $(E_{i}, \circ_{i}, e_{i})$ are FWL $F$-manifolds over  an $F$-manifold
$(M, *, e)$,  then  so is $(E_{1}\oplus E_{2}, \circ_{1}\oplus \circ_{2},
e_{1}\oplus e_{2})$.  If $(E_{i}, \circ_{i}, e_{i}, \nabla^{i})$ are FWL 
flat $F$-manifolds over  a flat $F$-manifold $(M, *, e,\nabla)$, then  so is 
$(E_{1}\oplus E_{2}, \circ_{1}\oplus \circ_{2},
e_{1}\oplus e_{2}, \nabla^{1}\oplus \nabla^{2})$, where 
the components $\nabla^{E}$ and $\mathcal S$ of $\circ_{1}\oplus \circ_{2}$ are 
given by $\nabla^{E_{1}}\oplus \nabla^{E_{2}}$ and  $\mathcal S_{1}\oplus \mathcal S_{2}.$
\end{lem}

We define  the {\cmssl generalized prolongation }
$(\mathbb{T}M = TM\oplus T^{*}M,  *^{\mathbb{T}}, e^{\mathbb{T}})$  
of  $(M, *, e, \nabla )$  as the direct sum 
of the tangent and cotangent prolongations.\\

{\bf Iterations of prolongations.}  Applying successively Proposition \ref{integr-comp}, we define  {\cmssl higher order tangent prolongations}  of an $F$-manifold.
The next lemma is essential in iterations involving the cotangent or generalized prolongations:
\begin{lem} Let $(M, *, e, \nabla )$ be a flat $F$-manifold. Then  $TM$, 
$T^{*}M$ and $\mathbb{T}M$ are  flat $F$-manifolds, with flat connections $\nabla^{T}$, 
$(\nabla^{T})^{*}$ and 
$\nabla^{T}\oplus (\nabla^{T})^{*}$.
\end{lem}

\begin{exa}{\rm  The $F$-manifold $T^{*}(TM)$ is the image of  the second tangent prolongation $T(TM)$ in the duality 
$I^{\nabla^{T}}$ of Theorem \ref{short} defined by the flat compatible connection $\nabla^{T}$ on $(TM, *^{T} , e^{T}).$
The $F$-manifold $T^{*}(T^{*}M)$ is the image of  the tangent prolongation $T(T^{*}M)$ of  the $F$-manifold $T^{*}M$ in 
the duality  $I^{(\nabla^{T})^{*}}$ of Theorem \ref{short} defined by the flat compatible connection  $(\nabla^{T})^{*}$ on 
$(T^{*}M, *^{T^{*}} , e^{T^{*}}).$}
\end{exa}

\begin{cor}\label{last}   i) Let $(M, *, e, \mathcal E)$ be an $F$-manifold with Euler field. Then all tangent 
prolongations  of $(M, *, e)$ inherit an Euler field.\

ii) Let $(M, *, e, \nabla  ,\mathcal E)$ be a flat  $F$-manifold with Euler field. Then 
all prolongations of $(M, *, e,\nabla )$ are flat $F$-manifolds with Euler fields. 
\end{cor}

Liana David: Institute of Mathematics  `Simion Stoilow' of the Romanian Academy,   Calea Grivitei no.\ 21,  Sector 1, 010702, Bucharest, Romania; 
liana.david@imar.ro
\end{document}